\documentclass[10pt]{amsart}
\usepackage{amssymb}
\usepackage{pst-node}
\usepackage{mathtools}
\usepackage{amsmath}
\usepackage{enumerate}
\usepackage{subcaption}
\usepackage{pgfplots}
\pgfplotsset{compat=1.18}
\usepgfplotslibrary{fillbetween}
\usepackage{tikz-cd}
\usepackage{float}
\usetikzlibrary{decorations.pathreplacing}
\usepackage{graphicx} 

\newcommand{\OmegaAD}{\Omega_{\alpha}^{\delta}}
\newcommand{\Ra}{R_{\alpha}}
\newcommand{\eps}{\varepsilon}
\newcommand{\J}{\mathcal{J}}
\newcommand{\Dioph}{\operatorname{Dioph}}
\newcommand{\dOm}{d_{\Omega^{\delta}_{\alpha}}}

\newtheorem{thm}{Theorem}[section]
\newtheorem{prop}[thm]{Proposition}
\newtheorem{lem}[thm]{Lemma}
\newtheorem{cor}[thm]{Corollary}

\theoremstyle{definition}
\newtheorem{definition}[thm]{Definition}

\newtheorem{rmk}[thm]{Remark}

\numberwithin{equation}{section}

\calclayout

\usepackage{hyperref}

\title{Lower Intermediate dimensions of classical Denjoy minimal sets}
\author{Milo Edwardes}
\date{\today}

\begin{document}

\begin{abstract}
    We study the $\theta$-intermediate dimension of the minimal set of the classical $C^{1+\delta}$ Denjoy map, an orientation-preserving diffeomorphism of the circle not conjugate to a rigid rotation. This follows on from work of Kra and Schmeling, who showed how the box and Hausdorff dimensions may differ, depending on the Diophantine properties of the rotation number. The $\theta$-intermediate dimension interpolates between these two dimensions. We use the well-understood gap structure of the irrational rotation map to develop a covering strategy to get upper bounds on the lower $\theta$-intermediate dimension, and use results of Banaji to show that these upper bounds are sharp for $\theta$ in an interval.
\end{abstract}
\maketitle
\section{Introduction}

In his work from 1881--1886, Poincaré \cite{HPoincare} gave a classification of all orientation-preserving homeomorphisms of $S^1=\mathbb{R}/\mathbb{Z}$, proving that any such homeomorphism $f$ has a unique associated real parameter $\alpha=\alpha(f)$ called its rotation number, and, if $\alpha$ is irrational, then $f$ is semiconjugate to the rigid rotation $\Ra:S^1\rightarrow S^1$ defined by $\Ra(x)=x+\alpha \bmod 1$. We recall that $f:X\rightarrow X$ is \textit{semiconjugate} to $g:Y\rightarrow Y$  if there exists some continuous surjective map $h:X\rightarrow Y$ such that $h \circ f=g\circ h$. We say $f$ and $g$ are \textit{conjugate} if $h$ is invertible. We say that $f$ is \textit{minimal} if the orbit of every $x\in X$ under $f$ is dense. For a homeomorphism $f$ we say that a subset $\Omega \subset X$ is \textit{minimal} if $\Omega $ is non-empty, compact, invariant and $f$ is minimal on $\Omega$, i.e. the orbit of every $x\in\Omega$ under $f$ is dense in $\Omega$. Poincaré's Theorem tells us that if $f$ is minimal, then it is conjugate to a rotation by some $\alpha$. Denjoy expanded on this \cite{ADenjoy} by showing that if $f$ has an irrational rotation number and the derivative of $f$ is of bounded variation, it is minimal and thus conjugate to a rotation. He also produced examples of $C^1$ diffeomorphisms not conjugate to rotations, details of which can be found in \cite{Melo1993}. For any $\delta<1$, Herman \cite{Herman1979} later gave examples of $C^{1+\delta}$ diffeomorphisms that are not conjugate to rotations.

We define a \textit{Denjoy map} as any orientation-preserving $C^1$ diffeomorphism with irrational rotation number not conjugate to a rotation. The examples of Denjoy maps presented by Denjoy and Herman, and in fact \emph{any} examples of Denjoy maps, have unique minimal sets which are Cantor sets (see also \cite{McDuff}). In 2002, Kra and Schmeling \cite{KraSchmeling} proved lower bounds for the box and Hausdorff dimensions of the minimal sets of $C^{1+\delta}$ Denjoy maps (Theorem \ref{thm:KS}): the bound for the box dimension depends only on $\delta$, and is attained by the classical Denjoy map (see also \cite{Norton1999}), whereas the bound for the Hausdorff dimension depends on the Diophantine class of the rotation number. See Section \ref{Sec:Prelims} for definitions of the classical Denjoy map and Diophantine class.

Falconer, Fraser and Kempton's paper \cite{Falconer2018} introduced the notion of $\theta$-intermediate dimension as a family of dimensions interpolating between the box and Hausdorff dimension, obtained by restricting to covers by sets with diameters in $[\Delta^{1/\theta},\Delta]$ for $\theta\in(0,1]$ (see Definition \ref{Def:IntermediateDim}). We refer to \cite{Falconer2021, FraserInterpolating, fraser2020assouaddimensionfractalgeometry} for overviews. The general theory has since been developed in several directions, including the attainable forms of the $\theta$-intermediate dimensions \cite{BanajiRutar}, generalised intermediate dimensions \cite{Banaji2023, BanajiRutarTroscheit}, projections and random images \cite{BurrellFalconerFraserProj, Feng, Burrell, FalconerFBM, FraserProjections, DawAttia} and measures \cite{AngeliniMolterTejada}. There has also been much work on the intermediate dimensions of dynamically defined sets, including Bedford--McMullen carpets \cite{Falconer2018, Kolossvary, BanajiKolossvary}, polynomial spirals \cite{BurrellFalconerFraserSpirals} and non-autonomous conformal sets \cite{MiaoWang}. Closest in spirit to this work is the result of Banaji and Fraser \cite{BanajiFraser}, who showed that the $\theta$-intermediate dimensions of the limit set of a conformal infinite iterated function system are the maximum of its Hausdorff dimension and the $\theta$-intermediate dimensions of the set of fixed points of the defining maps. Denjoy minimal sets form a natural family of a different kind, in which the disparity between the Hausdorff and box dimensions has an arithmetic origin: by \cite{KraSchmeling}, the box dimension of the minimal set depends only on $\delta$, whereas their lower bound for its Hausdorff dimension depends on the Diophantine class of $\alpha$. This paper examines this disparity. Our main result is the following.

\begin{thm}\label{Thm:main}
    Let $\delta\in(0,1)$ and $\alpha\in \Dioph(\nu)$ with  $\nu>1$. Let $\OmegaAD$ denote the minimal set of the $C^{1+\delta}$ classical Denjoy map with rotation number $\alpha$. Let $r_{\nu}$ denote the unique positive root of
    \begin{equation*}
        (1-\delta)r^2+\delta r-\nu=0.
    \end{equation*}
    Then, for $\theta\in \left[\frac{1}{r_{\nu}},1\right]$, the lower $\theta$-intermediate dimension of $\OmegaAD$ is
    \begin{equation*}
        \underline{\dim}_{\theta}\OmegaAD=
            \frac{\delta\theta}{\delta\theta+1-\delta}.
    \end{equation*}
    For $\theta\in\left(0,\frac{1}{r_{\nu}}\right)$, we have
    \begin{equation*}
        \max\left\{\dim_{\mathrm{H}}\OmegaAD,\frac{\delta\theta}{\delta\theta+1-\delta}\right\}\leq \underline{\dim}_{\theta}\OmegaAD\leq\frac{\delta r_{\nu}}{\nu}.
    \end{equation*}
    Moreover, the upper bound for $\theta\in\left(0,\frac{1}{r_{\nu}}\right)$ is optimal among single-level covers.
\end{thm}

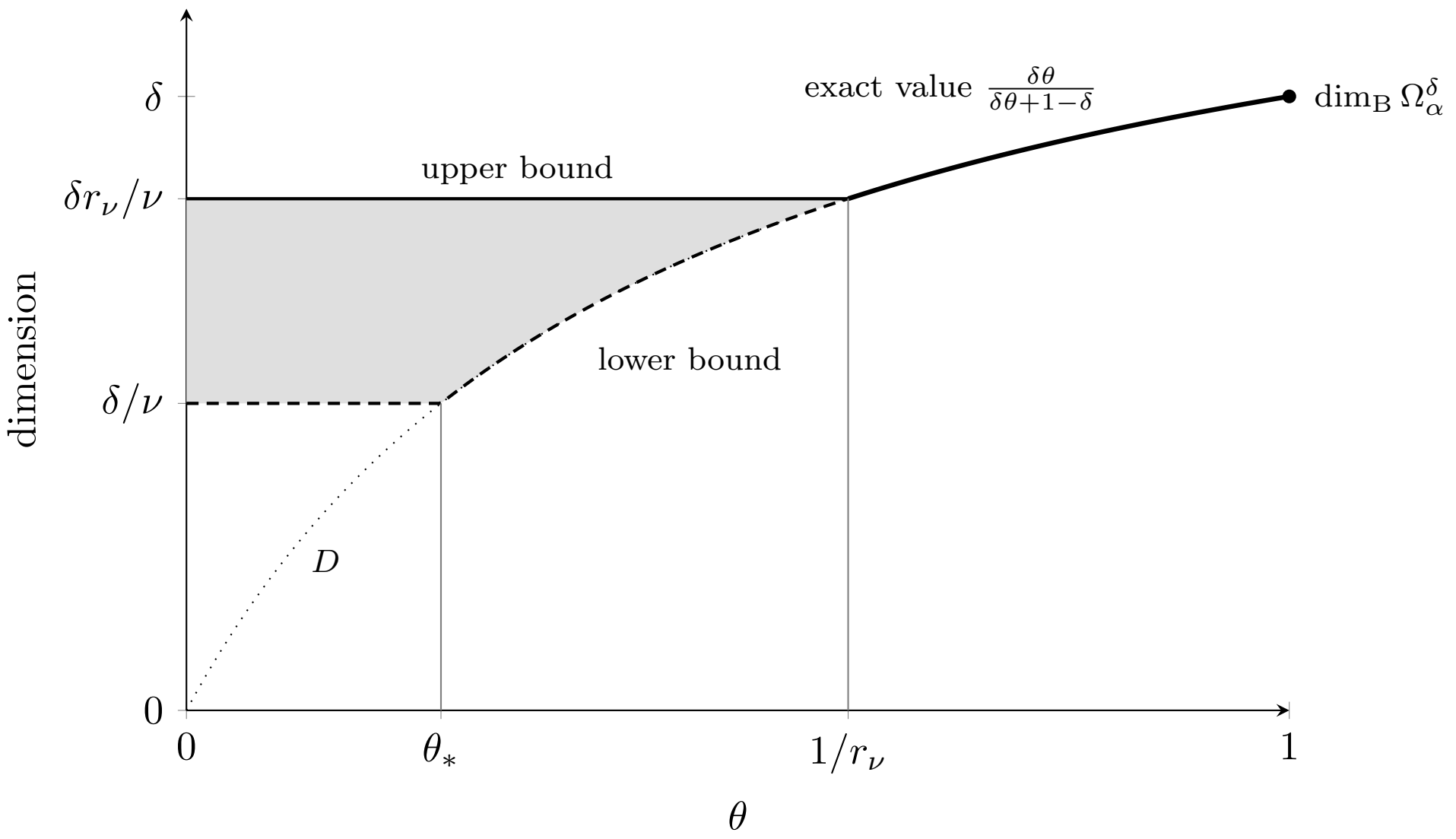
\begin{figure}[!htb]
\centering
\begin{tikzpicture}
\begin{axis}[
    width=0.8\linewidth, height=0.55\linewidth,
    xmin=0, xmax=1, ymin=0, ymax=0.8,
    xlabel={$\theta$}, ylabel={dimension},
    xtick={0,0.2308,0.6,1}, xticklabels={$0$,$\theta_{*}$,$1/r_{\nu}$,$1$},
    ytick={0,0.35,0.5833,0.7}, yticklabels={$0$,$\delta/\nu$,$\delta r_{\nu}/\nu$,$\delta$},
    axis lines=left, clip=false,
    every axis plot/.append style={samples=200},
]
\addplot[name path=upper, thick, domain=0:0.6] {7/12};
\addplot[name path=lower, thick, densely dashed, domain=0:0.6] {max(0.35, 0.7*x/(0.7*x+0.3))};
\addplot[gray!25] fill between[of=upper and lower];
\addplot[very thick, domain=0.6:1] {0.7*x/(0.7*x+0.3)};
\addplot[thin, dotted, domain=0:0.6] {0.7*x/(0.7*x+0.3)};
\addplot[only marks, mark=*, mark size=1.5pt] coordinates {(1,0.7)};
\addplot[gray, thin] coordinates {(0.6,0) (0.6,0.5833)};
\addplot[gray, thin] coordinates {(0.2308,0) (0.2308,0.35)};
\node[anchor=south,font=\footnotesize] at (axis cs:0.3,0.5833) {upper bound};
\node[anchor=north west,font=\footnotesize] at (axis cs:0.36,0.43) {lower bound};
\node[anchor=south east,font=\footnotesize] at (axis cs:0.84,0.665) {exact value $\frac{\delta\theta}{\delta\theta+1-\delta}$};
\node[anchor=west,font=\footnotesize] at (axis cs:0.1,0.17) {$D$};
\node[anchor=west,font=\footnotesize] at (axis cs:1.01,0.7) {$\dim_{\mathrm{B}}\OmegaAD$};
\end{axis}
\end{tikzpicture}
\caption{The bounds of Theorem \ref{Thm:main} on $\underline{\dim}_{\theta}\OmegaAD$ for $\delta=0.7$ and $\nu=2$, where $\theta_{*}=(1-\delta)/(\nu-\delta)$. The dotted curve is the value for the set $D$ of \cite{AngeliniMolter}.}
\label{Fig:IntDimBounds}
\end{figure}

We will be more explicit with what we mean by \emph{single-level covers} in Section \ref{Sec:UpperIntDimensions}. See Figure \ref{Fig:IntDimBounds} for a graph of this result.

Theorem \ref{Thm:main} describes how the lower intermediate dimension interpolates between the box and Hausdorff dimensions of $\OmegaAD$. For comparison, we recall the bounds of Kra and Schmeling on these dimensions.

\begin{thm}[\cite{KraSchmeling}]\label{thm:KS}
Let $\delta\in(0,1)$ and let $\alpha\in(0,1)$ be of Diophantine class $\nu$. Then the minimal set $\Omega$ of any orientation-preserving $C^{1+\delta}$ Denjoy map with rotation number $\alpha$ satisfies
\begin{equation*}
    \dim_{\mathrm{B}}\Omega\geq\delta \qquad\text{and}\qquad \dim_{\mathrm{H}}\Omega\geq\frac{\delta}{\nu}.
\end{equation*}
Moreover, $\dim_{\mathrm{B}}\OmegaAD=\delta$ for the classical Denjoy sequence.
\end{thm}

Angelini and Molter \cite{AngeliniMolter} recently studied sets closely related to Denjoy minimal sets.  Given a positive, non-increasing sequence $(a_{n})$ with $\sum_{n}a_{n}=1$, they considered the family of all closed subsets of $[0,1]$ whose complementary open intervals have lengths given by a rearrangement of $(a_{n})$, and showed that the smallest possible $\theta$-intermediate dimensions among them are attained by the countable set $D$ of tail sums of $(a_{n})$. See \cite{AngMolSurvey} for a detailed survey on this topic. When we introduce the construction of Denjoy minimal sets in Section \ref{SubSec:Construction} we will see that $\OmegaAD$ is a member of this family, and for the classical Denjoy sequence their results give $\underline{\dim}_{\theta}D=\delta\theta/(\delta\theta+1-\delta)$. In Theorem \ref{Thm:main}, we show using a result from \cite{Banaji2023} that
\begin{equation*}
    \underline{\dim}_{\theta}\OmegaAD\geq\frac{\delta\theta}{\delta\theta+1-\delta}, \qquad \theta\in(0,1],
\end{equation*}
in agreement with their result. This bound depends only on the lengths of the gaps, and not on the order in which the rotation $\Ra$ places them around the circle. For $\theta\geq 1/r_{\nu}$, our upper bound coincides with it, so that $\underline{\dim}_{\theta}\OmegaAD=\delta\theta/(\delta\theta+1-\delta)$, and at these scales $\OmegaAD$ behaves like the countable set $D$. For $\theta<1/r_{\nu}$, our upper bound is the constant $\delta r_{\nu}/\nu$, which depends on the Diophantine class $\nu$ of $\alpha$, and we show that no cover at a single level does better.

The paper is organised as follows. Section \ref{Sec:Prelims} contains the necessary preliminaries, introducing the construction of the classical Denjoy minimal set, as well as the required background from Diophantine approximation and fractal geometry which will be used throughout. Section \ref{Sec:Approx} employs the Three Gap Theorem, along with Diophantine information about the rotation number $\alpha$ in order to develop a systematic understanding of the minimal set at each level. This is then used to form a covering strategy used in the subsequent work. Finally, Section \ref{Sec:UpperIntDimensions} begins with some justification of our covering strategy before moving on to showing its optimality among single-level covers, and using it to deduce an upper bound on the lower $\theta$-intermediate dimension. We use a result from \cite{Banaji2023} to get lower bounds, showing that our upper bound is sharp for $\theta\geq\frac{1}{r_{\nu}}$. We end with a short remark on the upper $\theta$-intermediate dimension.

\subsection*{Acknowledgements} The author would like to thank Tom Kempton for introducing him to the problem, and Xiong Jin for helpful insights and support in the latter stages of completing this manuscript. Finally, the author thanks Jonathan Fraser, Nicolas Angelini and Ursula Molter for their helpful comments on an earlier version of this paper. This research was funded by the Heilbronn Institute for Mathematical Research (HIMR) through the Heilbronn Doctoral Partnership (HDP).

\subsection*{AI Statement} 
Claude (Anthropic) models have been used for literature search, proofreading, restructuring, and fixing errors and wording in earlier versions of this work. The entire manuscript was prepared and checked by the author.

\section{Preliminaries}\label{Sec:Prelims}
Throughout this work, we will use standard asymptotic notation. Given two real-valued functions on some domain $D$ $f,g:D\rightarrow\mathbb{R}$, we say that $f=O(g)$ if there exists some constant $C>0$ such that $f(x)\leq Cg(x)$ for all $x\in D$. We will say that $f=o(g)$ if $\frac{f(x)}{g(x)}\rightarrow 0$ as $x\rightarrow\infty$.  Similarly, for $D\subset \mathbb{R}$ some domain unbounded above, we will write $f=\Theta(g)$ if there exist constants $c,C>0$ and some $x_{0}\in D$ such that 
\begin{equation*}
    cg(x)\leq f(x) \leq Cg(x) \;\; \text{ for all }x\geq x_{0}.
\end{equation*}
\subsection{Construction of Denjoy maps}\label{SubSec:Construction}
As stated above, the fractals studied here are dynamically defined as the minimal set of a Denjoy map. To construct this, we introduce the following, as found in \cite{KraSchmeling}. For more details on the origins of this construction, see \cite{Melo1993,Katok1995,Navas}. This work will focus on Denjoy maps generated by classical Denjoy sequences.

\begin{definition}
    The \textit{classical Denjoy sequence}, denoted $D_{0}^{\delta}$, consists of lengths $$l_{n}=c_{\delta}(1+|n|)^{-1/\delta},$$ where $c_{\delta}^{-1}=\sum_{n\in\mathbb{Z}}(1+|n|)^{-1/\delta}$.
\end{definition}

Note that in the classical Denjoy sequence $l_{n}=l_{-n}$, and that the constant $c_{\delta}$ is a normalising constant, ensuring that $\sum_{n\in\mathbb{Z}} l_{n} = 1$. If instead we had that $\sum_{n\in\mathbb{Z}}l_{n}<1$, then our Denjoy minimal set would have positive Lebesgue measure and so questions of dimension would be uninteresting. The reason for this will become clear with the following description of the Denjoy construction. 

We are now ready to construct a classical Denjoy map. We note that we can generalise this construction for a broader family of Denjoy sequences, which we will not specify in this work (see \cite{KraSchmeling}). First, let $\alpha\in (0,1)$ be irrational, let $\delta\in (0,1)$ and let $(l_{n})$ be the classical Denjoy sequence of class $\delta$, $D_{0}^{\delta}$. To obtain a Denjoy map, we ``blow up" each point $n\alpha$ in the orbit of $0$ under $\Ra$ to an open interval of length $l_{n}\in D_{0}^{\delta}$ which we will call $J_{n}$ (and we will refer to them as \textit{Denjoy intervals}). At each step, we rescale the remaining part of the circle, maintaining a total length of 1.

In other words, we can think of this as a step-wise process in the following way. In the first step, we drop the interval $J_{0}$ of length $l_{0}$ over the point $0=\Ra^{0}(0)$, giving us a circle with total length $1+l_{0}$, before rescaling the original circle (not including the newly inserted interval) by $1-l_{0}$, leaving us with $J_{0}$ and an arc of length $1-l_{0}$. Next, we move to the point on $S^{1}$ where $\Ra(0)$ now lies, then place an interval $J_{1}$ of length $l_{1}$ at this point, so that the total length of the circle with this new interval is $1+l_{1}$. Finally, we rescale the complement of $J_{0}\cup J_{1}$ by the factor $(1-l_{0}-l_{1})/(1-l_{0})$, so that we have $J_{0}$, $J_{1}$ and a complementary arc of length $1-l_{0}-l_{1}$.

We continue in this manner, remembering to drop intervals on backward orbits of $\Ra$ as well as forward ones. After inserting the intervals $J_{n}$ with $|n|\leq N$, their complement has total length $1-\sum_{|n|\leq N}l_{n}$, which tends to $0$ as $N\to\infty$; in particular, the resulting minimal set has Lebesgue measure zero. Note that the rescaling of the complementary intervals at each step will change the relative position of iterates of $\Ra$ to $0$, i.e. we can imagine that we mark every point in the orbit of $\Ra(0)$ at the start of the process before adding in our intervals $J_{n}$ and rescaling the complementary intervals, which means that the rescaling will shift the iterates of $\Ra$ around the circle, but will preserve the cyclic ordering.

A Denjoy map $f$ is defined by mapping the interval $J_{n}$ diffeomorphically to $J_{n+1}$ with derivative 1 at the endpoints. By continuity, $f$ is defined on the entire circle and since $l_{n+1}/l_{n}\rightarrow1$, the maps sending $J_{n}\rightarrow J_{n+1}$ can be chosen with derivatives uniformly close to $1$, and $f$ can be chosen to be $C^{1+\delta}$ (see \cite{Norton1999,KraSchmeling,Herman1979}). Since the ordering of iterates of $\Ra$ determines the ordering of the intervals $J_{n}$, we have that $f$ will have rotation number $\alpha$ by construction. We define the \emph{classical Denjoy minimal set} of $f$ (depending on $D_{0}^{\delta}$) as follows:
\begin{equation*}
    \OmegaAD(D_{0}^{\delta})=\OmegaAD = S^1 \setminus \bigcup_{n\in\mathbb{Z}} J_{n}.
\end{equation*}
This is the unique minimal set of $f$ and is a Cantor set. We will also define a truncated version of this set as
\begin{equation*}
    \mathcal{J}_{N} = S^{1}\setminus\bigcup_{-N\leq n\leq N} J_{n}.
\end{equation*}

A defining feature of a Denjoy map is the semiconjugacy relating it to its associated irrational rotation. For a Denjoy map $f:S^1 \rightarrow S^1$, we have that $f$ is semiconjugate to $\Ra:S^1\rightarrow S^1$ via $h:S^1 \rightarrow S^1$ satisfying $h\circ f=\Ra \circ h$ by Poincar\'e's Theorem. Its restriction $h|_{\OmegaAD}$ is one-to-one except over the orbit of $0$, where each fibre consists of the two endpoints of $J_{n}$. In other words, we describe the fibres of $h$ as, for $x\in \{n\alpha \bmod 1 :n\in\mathbb{Z}\}$, there exists some $k\in\mathbb{Z}$ such that $x=k\alpha\bmod 1$ and we define the preimage $h^{-1}(x)$ to be the doubleton set of endpoints of the interval $J_{k}$. We label this pair of preimages as $h^{-1}(x)=\{g_{r}(x),g_{l}(x)\}$ as the right and left preimages respectively. For $x\in S^1\setminus \{n\alpha \bmod 1 :n\in\mathbb{Z}\}$, $h^{-1}(x)$ is well defined, in other words, $h^{-1}(x)=g_{r}(x)=g_{l}(x)$. We define $(x,y)$ to be the shorter oriented arc between $x$ and $y$, taking the positively oriented arc starting at $x$ in the antipodal case. Then, for a positively oriented arc, we define a metric

\begin{equation*}
   \dOm(g_{r}(x),g_{l}(y))= \sum_{\{n:n\alpha\in (x,y)\}} |J_{n}| + \bigg(1-\sum_{n\in\mathbb{Z}} |J_{n}|\bigg)\|x-y\|
\end{equation*}
where the second term vanishes since $\sum_{n}l_{n}=1$. For a negatively oriented arc $(x,y)$, we analogously define $d_{\Omega^{\delta}_{\alpha}}(g_{l}(x),g_{r}(y))=\sum_{\{n:\,n\alpha\in(y,x)\}}l_{n}$.

In other words, $\dOm$ is the arc-length metric of the circle obtained by inserting the intervals $J_{n}$, restricted to $\OmegaAD$. In particular, for an arc $(x,y)$ of $S^{1}$, the corresponding interval of $\OmegaAD$ runs between the facing endpoints $g_{r}(x)$ and $g_{l}(y)$, and its length $\sum_{n\alpha\in(x,y)}l_{n}$ does not include the Denjoy intervals at $x$ and $y$.

We give the following figure as an example.

\begin{figure}[!ht]
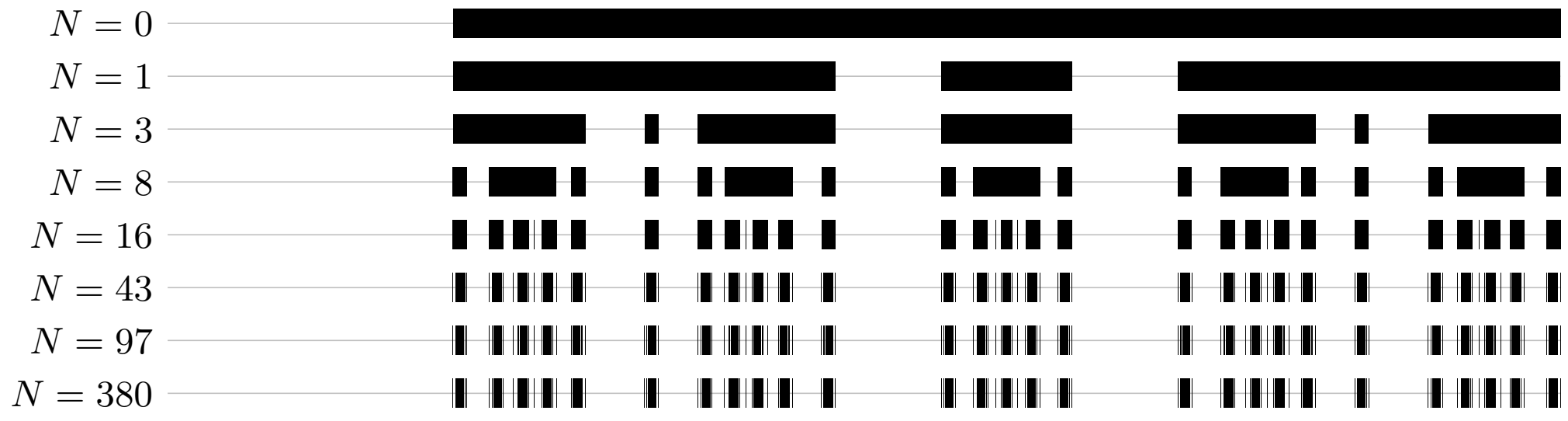

\centering

\caption{The sets $\J_{N}$ (black) for $\delta=0.7$ and $\alpha=[0;2,2,5,27,734,\ldots]\in\Dioph(2)$, drawn in the $\dOm$ metric.}
\label{Fig:DenjoyLevels}
\end{figure}

\subsection{Dimension}
   The results of this paper are based on the dimension properties of the Denjoy minimal set. As such, we define the notions of dimension which will come up in the following work. Throughout, $F$ denotes a non-empty bounded subset of a metric space $(X,d)$. We will apply these definitions with $X=S^{1}$, equipped with the arc-length metric of the circle constructed in Section \ref{SubSec:Construction}, and $F=\OmegaAD$. We will write $N_{F}(\eps)$ for the smallest number of sets of diameter $\eps$ needed to cover $F$.
\begin{definition}
    For a bounded subset $F$ of a metric space $X$, the \textit{upper box dimension} of $F$, denoted $\overline{\dim}_{\mathrm{B}}F$ is given by
    \begin{equation*}
        \overline{\dim}_{\mathrm{B}}F:= \limsup_{\eps\rightarrow 0 }\frac{\log N_{F}(\eps)}{\log 1/\eps}.
    \end{equation*}
    We define the \textit{lower box dimension} similarly, denoted $\underline{\dim}_{\mathrm{B}}F$, as
    \begin{equation*}
        \underline{\dim}_{\mathrm{B}}F:= \liminf_{\eps\rightarrow 0 }\frac{\log N_{F}(\eps)}{\log 1/\eps}.
    \end{equation*}
\end{definition}

The next notion of dimension we introduce is the Hausdorff dimension. To do this, we first introduce the Hausdorff measure. For some subset $U$ of $X$, we let diam$(U):= \sup \{d(x,y):x,y\in U\}$ denote the diameter of $U$.

\begin{definition}
    Let $s\in[0,\infty)$. The $s$\textit{-dimensional Hausdorff measure} $\mathcal{H}^{s}(F)$ of a subset $F$ of a metric space $X$ is defined by
    \begin{align*}
        \mathcal{H}_{\eps}^{s}(F)=\inf\left\{\sum_{i}(\operatorname{diam}(U_{i}))^{s}:  F\subset \bigcup_{i}U_{i} \text{ a countable cover,  and } \sup_{i}\operatorname{diam}(U_{i})<\eps\right\},
    \end{align*}
    and 
    \begin{equation*}
        \mathcal{H}^{s}(F):=\lim_{\eps\downarrow 0}\mathcal{H}^{s}_{\eps}(F).
    \end{equation*}
\end{definition}

We note that there is a unique critical value $s_{0}=s_{0}(F)$ such that

\begin{equation*}
    \mathcal{H}^{s}(F)=
    \begin{cases}
        \infty &\text{for }s<s_{0} \\
        0  &\text{for } s>s_{0}. 
    \end{cases}
\end{equation*}
\begin{definition}
    The unique critical value $s_{0}$ as described above is defined to be the Hausdorff dimension of $F$, denoted $\dim_{\mathrm{H}}F$.
\end{definition}

We now define what we mean by $\theta$-intermediate dimensions. As mentioned in the introduction, Falconer, Fraser and Kempton first described this in \cite{Falconer2018}. In the following definition, we write $|U|:=\operatorname{diam}(U)$.

\begin{definition}\label{Def:IntermediateDim}
      For $\theta \in(0,1] $, the \textit{upper/lower $\theta$-intermediate dimensions} of a bounded $F\subset X$ are as follows:
\begin{multline*}
    \overline{\dim}_{\theta}F:=\inf\bigg\{s\geq 0 :  \forall \; \varepsilon>0 \; \exists \; \Delta_{0}\in (0,1] \text{ such that } \forall \; \Delta\in (0,\Delta_{0})\; \exists  \text{ a cover } \{U_{i}\}\\ \text{ of } F \text{ such that } 
      \Delta^{1/\theta}\leq |U_{i}|\leq \Delta\; \forall \;i, \text{ and } \sum_{i} |U_{i}|^{s}\leq \varepsilon\Bigg\}.
\end{multline*}
\begin{multline*}
    \underline{\dim}_{\theta}F:=\inf\Big\{s\geq 0:\forall\; \varepsilon>0\; \forall\;\Delta_{0}>0,\; \exists\;\Delta\in (0,\Delta_{0})\text{ and a cover }\{U_{i}\} \\\text{ of } F \text{ such that } \Delta^{1/\theta}\leq |U_{i}|\leq \Delta \;\forall \; i, \text{ and } \sum_{i}|U_{i}|^{s}\leq\varepsilon \Big\}.
\end{multline*}
We define $\overline{\dim}_{0 }F=\underline{\dim}_{0}F:=\dim_{\mathrm{H}}F$.
\end{definition}

When $F=\OmegaAD$, we write $|U|_{\dOm}$ for the diameter, and reserve $|I|$ for the Lebesgue length of an arc $I\subset S^1$.

The $\theta$-intermediate dimensions bridge the gap between box and Hausdorff dimensions in the following way: box dimension requires us to ``measure'' our fractal at a given level by covering sets of a fixed diameter $r$. Hausdorff dimension views the fractal at the same level but allows us to use covering sets $U$ such that $0<|U|\leq r$. $\theta$-intermediate dimensions provide a compromise between these cases, saying that our covering sets $U_{i}$ may have diameters ranging from $r^{1/\theta}\leq |U_{i}|\leq r$. It clearly follows that $\overline{\dim}_{1}F=\overline{\dim}_{\mathrm{B}}F$ and $\underline{\dim}_{1}F=\underline{\dim}_{\mathrm{B}}F$ for any bounded $F$. It is also simple to see that for $\theta<\theta'$ we get that $\overline{\dim}_{\theta}F\leq\overline{\dim}_{\theta'}F$ and $\underline{\dim}_{\theta}F\leq\underline{\dim}_{\theta'}F$, in other words, the maps $\theta\mapsto \overline{\dim}_{\theta}F$ and $\theta\mapsto \underline{\dim}_{\theta}F$ are non-decreasing in $\theta\in(0,1]$. We refer the reader to \cite{Banaji2023,FalconerFracGeomFound,Falconer2021,fraser2020assouaddimensionfractalgeometry} for a more detailed study on fractal dimension.

\subsection{Diophantine Approximation}
    We now introduce the basic theory and notation of Diophantine approximation and numbers. The results and definitions in this subsection can be found in \cite{KhinchinCF,Harrison1989,Lang1966,KraSchmeling,vanRavenstein_1988}. Refer to these citations for more details.
    
    For every irrational number $\alpha$, we can express $\alpha $ as an \textit{infinite continued fraction} as follows:
    \begin{equation*}
        \alpha=[a_{0}; a_{1}, a_{2}, a_{3},\dots] = a_{0} + \cfrac{1}{a_{1}+\cfrac{1}{a_{2}+\cfrac{1}{a_{3}+\ddots}}}
    \end{equation*}
    where $a_{i}\in\mathbb{N}$ for all $i\geq1$, and $a_{0}\in\mathbb{Z}$. In this work, we will restrict to $\alpha\in(0,1)$ and so can assert that $a_0=0$ and $a_{n}\in\mathbb{N}$ for all $n\geq 1$.
    We define the \textit{convergents} of $\alpha$ to be the sequence of irreducible fractions
    \begin{equation*}
        \frac{p_{n}}{q_{n}}=[a_{0}; a_{1}, a_{2},\dots, a_{n}]=a_{0}+\cfrac{1}{a_{1}+\cfrac{1}{a_{2}+\cfrac{1}{\ddots \cfrac{1}{a_{n-1}+\cfrac{1}{a_{n}}}} }}.
    \end{equation*}
We can define convergents recursively, for $n\geq 2$, as follows:
\begin{align*}
    p_{0}&=0, p_{1} =1; &&\text{ and } \;\;\;\;\;\;\;p_{n}= a_{n}p_{n-1}+p_{n-2}
    \\ q_{0}&=1,q_{1}=a_{1}; &&\text{ and }\;\;\;\;\;\;\; q_{n}= a_{n}q_{n-1}+q_{n-2}.
\end{align*}
This way of expressing continued fraction convergents gives rise to the following definition.
\begin{definition}\label{partialconvergents}
    The \textit{intermediate convergents} of an irrational number $\alpha$ are the irreducible fractions
    \begin{equation*}
        \frac{p_{n,i}}{q_{n,i}}=\frac{p_{n-2}+ip_{n-1}}{q_{n-2}+iq_{n-1}} = [a_{0}; a_{1}, a_{2},\cdots, a_{n-1},i].
    \end{equation*}
\end{definition}
In particular, we will define $p_{n,i}=ip_{n-1}+p_{n-2}$ and $q_{n,i}=iq_{n-1}+q_{n-2}$ for $0\leq i\leq a_{n}$. From this definition, clearly $q_{n,0}=q_{n-2}$ $q_{n,a_{n}}=q_{n}$.
We let $\|\cdot\|$ denote the distance to the nearest integer. We state the following result about this norm on our continued fraction convergents. Later on, we will want to understand quantities of the form $\|q_{n}\alpha\|$.
A commonly known fact about $\|q_{n}\alpha\|$ is that
\begin{equation}\label{eqn:CFnormineqality}
    \frac{1}{q_{n}(a_{n+1}+2)}\leq \|q_{n}\alpha\|< \frac{1}{a_{n+1}q_{n}}.
\end{equation}

We also have (see \cite{Lang1966}), that
\begin{equation}\label{eqn:CFnorminequality2}
    \frac{1}{2q_{n+1}}< \|q_{n}\alpha\|<\frac{1}{q_{n+1}}.
\end{equation}

\begin{prop}\label{2GapCase}
    For $n\in\mathbb{N}$, we have that for $\alpha\in(0,1)\setminus\mathbb{Q}$ with sequence of convergent denominators $\{q_{n}\}_{n\in\mathbb{N}}$
    \begin{equation}\label{2GapCaseEqn}
        q_{n}\|q_{n+1}\alpha\| + q_{n+1}\|q_{n}\alpha\| = 1
    \end{equation}
\end{prop}
Harrison gives a proof of this result in \cite{Harrison1989}.

Next, we define a property of an irrational number that will prove central to the dimension results given in what follows.

\begin{definition}
    An irrational number $\alpha$ is of \textit{Diophantine class} $\nu=\nu(\alpha)\in\mathbb{R}^+$ if 
    \begin{equation}
        \|q\alpha\|<\frac{1}{q^{\mu}}
    \end{equation}
    has infinitely many solutions for $q\in\mathbb{N}$ for $\mu<\nu$ and at most finitely many for $\mu>\nu$. For $\alpha$ of Diophantine class $\nu$, we write $\alpha\in \Dioph(\nu)$.
\end{definition}
The Diophantine class of an irrational number measures the strength of the exceptional rational approximations to $\alpha$. There is a quantity closely related to the Diophantine class called the \textit{irrationality measure} or \textit{irrationality exponent} defined using a more classical formulation from rational approximation, as can be found in \cite{Borwein1987}. A famous result of Dirichlet's tells us that all irrational numbers are of Diophantine class at least $1$ (and so of irrationality measure at least $2$). In this paper, we will primarily be interested in $\alpha\in\Dioph(\nu)$ with $\nu>1$. When $\nu=1$, Theorem \ref{thm:KS} gives $\dim_{\mathrm{H}}\OmegaAD=\dim_{\mathrm{B}}\OmegaAD=\delta$, so all of the $\theta$-intermediate dimensions of $\OmegaAD$ equal $\delta$.

The final thing to remark on in this subsection is a property of the Diophantine class which will be useful to us in our dimension arguments. A defining property of Diophantine class presented by Sondow in \cite{sondow2004irrationalitymeasuresirrationalitybases} is that for $\alpha\in \Dioph(\nu)$ we have that 

\begin{equation}
    \limsup_{n\rightarrow\infty} \frac{\log q_{n+1}}{\log q_{n} } = \nu.
\end{equation}
 This tells us that there are subsequences along which $q_{n+1}=q_{n}^{\nu+o(1)}$. In particular, using \eqref{eqn:CFnorminequality2}, we can show that there exists some subsequence $\{n_{k}\}_{k\in\mathbb{N}}\subset \mathbb{N}$ such that for any 
$\varepsilon>0$, there exists $K_{\eps}\in \mathbb{N}$ such that for all $k\geq K_{\eps}$, 
\begin{equation}\label{diophclassapprox}
   \frac{1}{2}q_{n_{k}}^{\nu -\eps}\leq q_{n_{k}+1} \leq q_{n_{k}}^{\nu+\varepsilon}.
\end{equation}
We will use a version of this later, which we state in the following lemma.
\begin{lem}\label{lem:DiophClass}
    Let $\alpha\in \Dioph(\nu)$, $\nu>1$, be irrational with $\{q_{n}\}$ its sequence of convergent denominators. Then, for every $\eps>0$ there exists $n_{\eps}\in \mathbb{N} $ such that for all $n\geq n_{\eps}$
    \begin{equation}
        q_{n+1} < q_{n}^{\nu+\eps}.
    \end{equation}
\end{lem}
\begin{proof}
    Since $\alpha\in \Dioph(\nu)$, we have that for every $\eps>0$, $$\|q\alpha\|<\frac{1}{q^{\nu+\eps}}$$ has at most finitely many solutions for integers $q$. This means that for all but finitely many $n\in\mathbb{N}$, $\|q_{n}\alpha\|>\frac{1}{q_{n}^{\nu+\eps}}$, and in particular, this means that for each $\eps>0$ there exists some $n_{\eps}\in\mathbb{N}$ such that for all $n\geq n_{\eps}$, $\|q_{n}\alpha\|>\frac{1}{q_{n}^{\nu+\eps}}$. We also  have by \eqref{eqn:CFnorminequality2} that $\|q_{n}\alpha\|<\frac{1}{q_{n+1}}$, so for $n\geq n_{\eps}$ we have that
    \begin{equation*}
        \frac{1}{q_{n}^{\nu+\eps}}<\|q_{n}\alpha\|<\frac{1}{q_{n+1}}
    \end{equation*}
    giving us that $q_{n+1}<q_{n}^{\nu+\eps}$ for all $n\geq n_{\eps}$.
\end{proof}
\subsection{Three Gap Theorem}
Throughout, we let $\{q_{n}\}_{n\in\mathbb{N}}$ define the sequence of denominators of the convergents of the continued fraction expansion of some (specified) irrational number $\alpha$. Many of the arguments in this paper will use the gap structure of the circle as partitioned by the orbit of $0$ under an irrational rotation. To use it later, we begin by stating the Three Gap Theorem as conjectured by Steinhaus, found in \cite{vanRavenstein_1988} by van Ravenstein. For a modern reference on this topic, see \cite{MarklofStrombergsson}.

For $x\in\mathbb{R}$, let $\{x\}$ denote the fractional part of $x$, and $\lfloor x \rfloor$ the integer part of $x$, so that $x=\lfloor x \rfloor+\{x\}$. 
Consider the sequence $\omega_{N}(\alpha):= \{ \{n\alpha\} = n\alpha \bmod 1\ : n= 0,1,2,...,N-1\}$. The Three Gap Theorem states that the unit interval (or $\mathbb{R}/\mathbb{Z}$, or $S^{1}$) is partitioned into arcs or gaps of at most three, and at least two, different lengths for any irrational $\alpha$ and any $N\in\mathbb{N}$. When visualising this problem, it is more helpful to consider the elements of $\omega_{N}(\alpha)$ in ascending order of magnitude. So let $(\{u_{j}(N)\alpha\})$ for $j=1,2,...,N$ be the ordered sequence of points in $\omega_{N}(\alpha)$. This gives us the sequence $\{u_{1}(N),u_{2}(N),...,u_{N}(N)\}$ of natural numbers $\{0,1,...,N-1\}$ such that $\{u_{j}(N)\alpha\}<\{u_{j+1}(N)\alpha\}$ for each $j$. In later arguments, we will abbreviate $u_{j}(N)$ to write $u_{j}$, where $N$ will have been fixed. As promised, fix $N\in\mathbb{N}$ and let $u_{j}$ denote the natural number representing the $j^{\text{th}}$ appearing iterate of $\Ra$ for some irrational $\alpha$. We define operators $Pre$ and $Suc$ such that $u_{j}=Pre(u_{j+1}) = Suc(u_{j-1})$, for $2\leq j \leq N$. For the endpoints, we define $Suc(u_{N})=u_{1}$ and $Pre(u_{1})=u_{N}$.

The Steinhaus conjecture, commonly known as the Three Gap Theorem, is stated as follows.

\begin{thm}[Three Gap Theorem \cite{vanRavenstein_1988}] \label{SteinhausConj}
    \begin{equation*}
        Suc(m)-m =
        \begin{cases}
            u_{2}, & 0\leq m < N-u_{2}, \\
            u_{2}-u_{N}, & N-u_{2} \leq m< u_{N}, \\
            -u_{N}, & u_{N} \leq m <N.
        \end{cases}
    \end{equation*}
\end{thm}

This tells us that the difference between the natural numbers representing successive points in $\{\Ra^{i}(0)\}_{i=0}^{N-1}$ as they appear in $S^1$ can take at most one of three values. While this information is useful, we are interested in the lengths of these gap intervals. As such, we state the following corollary providing this information. We denote the length of an interval in $S^1$ by $|\cdot|$ and let the distance from the nearest integer norm be given by $\|\cdot\|$. 

\begin{cor}
   \begin{equation*}
        \left|\left(\{m\alpha\},\{Suc(m)\alpha\}\right)\right| =
        \begin{cases}
            \|u_{2}\alpha\|, & 0\leq m < N-u_{2}, \\
            \|(u_{2}-u_{N})\alpha\|, & N-u_{2} \leq m< u_{N}, \\
            \|u_{N}\alpha\|, & u_{N} \leq m <N.
        \end{cases}
    \end{equation*}

\end{cor}

The next theorem indicates what these values are, and under which circumstances.  

\begin{thm}[T. van Ravenstein \cite{vanRavenstein_1988}]\label{Thm:vanRavenstein}
    For $q_{n,i-1} < N \leq q_{n,i}$, $2\leq i \leq a_{n}$ (and $n\geq 2$), we have that 
    $$ u_{2} =
    \begin{cases}
        q_{n-1}, & n \text{ odd,} \\
        q_{n,i-1}, & n \text{ even},
    \end{cases}
    \;\;\;
     u_{N} =
    \begin{cases}
        q_{n,i-1}, & n \text{ odd,} \\
        q_{n-1}, & n\text{ even}.
    \end{cases}
    $$
    For $q_{n-1}<N\leq q_{n,1}$ ($n\geq 2)$,
    $$
    u_{2} =
    \begin{cases}
        q_{n-1}, & n \text{ odd,} \\
        q_{n-2}, & n \text{ even,}       
    \end{cases}
    \;\;\;
    u_{N}= 
    \begin{cases}
        q_{n-2}, & n \text{ odd,} \\
        q_{n-1}, & n \text{ even.}
    \end{cases}
    $$
    For $N\leq q_{1}$ we have that $u_{j} = j-1$ for all $j= 1,2,...,N$.
\end{thm}

In our discussion, we will focus on the case where $n$ is odd since the even case is symmetric. So, putting these two results together, we have the following.

\begin{cor}\label{3gaplengths}
    For $q_{n-1}<N\leq q_{n,1}$, $(n\geq 2)$, $n$ odd, we have that
     \begin{equation*}
        \left|\left(\{m\alpha\},\{Suc(m)\alpha\}\right)\right| =
        \begin{cases}
            \|q_{n-1}\alpha\|, & 0\leq m < N-q_{n-1}, \\
            \|(q_{n-1}-q_{n-2})\alpha\|, & N-q_{n-1} \leq m< q_{n-2}, \\
            \|q_{n-2}\alpha\|, & q_{n-2} \leq m <N.
        \end{cases}
    \end{equation*}
    For $q_{n,i-1} < N \leq q_{n,i}$, $2\leq i \leq a_{n}$ (and $n\geq 2$ and odd), we have that
     \begin{equation*}
        \left|\left(\{m\alpha\},\{Suc(m)\alpha\}\right)\right| =
        \begin{cases}
            \|q_{n-1}\alpha\|, & 0\leq m < N-q_{n-1}, \\
            \|(q_{n-1}-q_{n,i-1})\alpha\|, & N-q_{n-1} \leq m< q_{n,i-1}, \\
            \|q_{n,i-1}\alpha\|, & q_{n,i-1} \leq m <N.
        \end{cases}
    \end{equation*}
\end{cor}

The results given above tell us about the gap structure between forward orbits, i.e. elements of the set $\{\Ra^{i}(0)\}_{i=0}^{N-1}\subset S^1$. The purpose of understanding this gap structure is to transfer this understanding to make sense of the interval structure of $\mathcal{J}_{n}$ - the Denjoy minimal set at some finite level. In constructing the Denjoy minimal set, we drop intervals whose size is indicated by our classical Denjoy sequence according to the forward and backward orbits of $\Ra$. As a result of this, the interval structure of $\mathcal{J}_{N}$ resembles the gap structure of $\{\Ra^{i}(0)\}_{i=-N}^{N}$. As argued in \cite{UrbanskiPawalek}, the Three Gap Theorem suffices in describing this gap structure by simply rotating $\Ra^{-N}(0)$ until it is at $0$. While the sizes of these intervals are preserved by this rigid rotation, the hitting times of intervals by iterates of $\Ra$ will not be. This property will be central to how we measure intervals with our metric on the minimal set, so restating the gap structure given in Corollary \ref{3gaplengths} where we have both forward and backward orbits, is useful. We note that for $N=2K+1$ where $N,K\in\mathbb{N}$, the gap structure of $\{\Ra^{i}(0)\}_{i=0}^{2K}$ is the same as the gap structure of $\{\Ra^{i}(0)\}_{i=-K}^{K}$, and so we rewrite Corollary \ref{3gaplengths} as follows:

\begin{cor}\label{Cor:Integer3GapThrm}
    Consider the set of points $\{\Ra^{i}(0)\}_{i=-K}^{K} \subset S^{1}$ partitioning the circle. Then for $q_{n-1}<2K+1\leq q_{n,1}$ with $n\geq 2$ and odd, we have that
    \begin{equation*}
        \left|\left(\{m\alpha\},\{Suc(m)\alpha\}\right)\right| =
        \begin{cases}
            \|q_{n-1}\alpha\|, & -K\leq m < K-q_{n-1}+1, \\
            \|(q_{n-1}-q_{n-2})\alpha\|, & K-q_{n-1} +1 \leq m< q_{n-2} -K, \\
            \|q_{n-2}\alpha\|, & q_{n-2} -K \leq m <K+1.
        \end{cases}
    \end{equation*}
    For $q_{n,i-1}<2K+1\leq q_{n,i}$, $2\leq i \leq a_{n}$ ($n\geq 2$ and odd), we have that 
    \begin{equation*}
        \left|\left(\{m\alpha\},\{Suc(m)\alpha\}\right)\right| =
        \begin{cases}
            \|q_{n-1}\alpha\|, & -K\leq m < K-q_{n-1}+1, \\
            \|(q_{n-1}-q_{n,i-1})\alpha\|, & K-q_{n-1}+1 \leq m< q_{n,i-1}-K, \\
            \|q_{n,i-1}\alpha\|, & q_{n,i-1}-K \leq m <K+1.
        \end{cases}
    \end{equation*}
\end{cor}

\section{Approximations on interval sizes}\label{Sec:Approx}
In this section, we estimate the $\dOm$-lengths of each complementary interval of $\J_{N}$ in terms of its first hitting time by the orbit of $\Ra(0)$ (Proposition \ref{Prop:BigAppendixProp}), and deduce that intervals of the same type have comparable $\dOm$-lengths (Lemma \ref{BigLem}).

\subsection{Interval types and hitting times}
The Three Gap Theorem gives us the gap structure of the circle as partitioned by the irrational rotation map. The purpose of this section is to translate the information about the gap structure in $S^1$ to the interval structure of our minimal set $\OmegaAD$. We will show that the size of an interval in the minimal set will correspond to some function of its preimage under the semiconjugacy relating the Denjoy map to its associated irrational rotation.

\begin{definition}
    Let $N\in\mathbb{N}$. Each complementary interval of $\J_{N}$ corresponds, via the semiconjugacy $h:\OmegaAD\rightarrow S^1$, to a gap of $S^1$ as partitioned by $\{\Ra^{r}(0)\}_{|r|\leq N}$. We call the complementary intervals of $\J_{N}$ corresponding to gaps of length $\|x\alpha\|$ the \emph{$x$-type} intervals, and we let $\mathcal{I}_{x}$ denote the set of $x$-type intervals of $\J_{N}$. Where convenient, we identify an $x$-type interval with the corresponding arc of $S^1$.
\end{definition}

Referring back to Corollary \ref{Cor:Integer3GapThrm}, we see that $x$ will be some continued fraction convergent or intermediate convergent denominator. In what follows, we let $|\cdot|_{\dOm}$ denote the length attributed to an interval by the $\dOm$ metric when taking the relevant right and left preimages of the interval's endpoints under the semiconjugacy $h$. For an open arc $I\subset S^1$, we write
\begin{equation*}
    H(I):=\{r\in\mathbb{Z}:r\alpha\in I\} \;\;\; \text{ and } \;\;\; \rho(I):=\min\{|r|:r\in H(I)\}
\end{equation*}
for the set of hitting times of $I$, and its first hitting time, respectively. The arguments below use work found in \cite{vanRavenstein_1988} to estimate hitting times.

\subsection{Hitting times and \texorpdfstring{$\dOm$}{dOm}-lengths}
\begin{lem}\label{Lem:HittingSeparation}
Let $\alpha\in(0,1)$ be irrational, and let $I\subset S^{1}$ be an open arc.
\begin{enumerate}[(i)]
    \item Let $m\geq1$. If $|I|\leq\|q_{m-1}\alpha\|$, then any two distinct $r,r'\in H(I)$ satisfy $|r-r'|\geq q_{m}$. In other words, every set of $q_{m}$ consecutive integers contains at most one element of $H(I)$.
    \item Let $m\geq2$. If $|I|\leq\|q_{m-1}\alpha\|+\|q_{m}\alpha\|$, then every set of $q_{m}$ consecutive integers contains at most two elements of $H(I)$.
\end{enumerate}
\end{lem}
 
\begin{proof}
For (i), since $r\alpha$ and $r'\alpha$ lie in the open arc $I$, whose length is at most $\|q_{m-1}\alpha\|\leq\frac{1}{2}$, we have $\|(r-r')\alpha\|<|I|\leq\|q_{m-1}\alpha\|$. By the best approximation property of the convergents of $\alpha$, the next integer $k$ for which $\|k\alpha\|<\|q_{m-1}\alpha\|$ is $q_{m}$, so we have that $|r-r'|\geq q_{m}$.
 
For (ii), let $x\in\mathbb{Z}$. The points $r\alpha$ with $x\leq r<x+q_{m}$ are a rotation of the points $r\alpha$ with $0\leq r<q_{m}$, and so, by Theorem \ref{Thm:vanRavenstein}, they partition $S^{1}$ into gaps of lengths $\|q_{m-1}\alpha\|$ and $\|q_{m-1}\alpha\|+\|q_{m}\alpha\|$. Any three of these points therefore span an arc of length at least $2\|q_{m-1}\alpha\|>\|q_{m-1}\alpha\|+\|q_{m}\alpha\|\geq|I|$, so the open arc $I$ contains at most two of them.
\end{proof}
 The next lemma collects the bounds on first hitting times used in the rest of this section.

\begin{lem}\label{Lem:BoundsOnRho}
    Let $I\subset S^1$ be an open arc. Then
    \begin{enumerate}[(i)]
        \item If $I$ is a complementary interval of $\J_{N}$, then $\rho(I)>N$.
        \item If, for some $L\in\mathbb{N}$, the points $\{\Ra^{r}(0)\}_{r=0}^{L}$ partition $S^1$ into gaps shorter than $|I|$, then $\rho(I)\leq \lceil \frac{L}{2}\rceil$.
        \item If $I$ is a complementary interval of $\J_{N}$ and $|I|\leq \|q_{m}\alpha\|$ for some $m\geq 0$, then $\rho(I)\geq q_{m+1}-N$.
    \end{enumerate}
\end{lem}
\begin{proof}
    Firstly, we note that $\rho(I)>N$ since $I$ contains no points $r\alpha$ with $|r|\leq N$, proving (i). 
    
    If we have that for $L\in\mathbb{N}$, $\{\Ra^{r}(0)\}_{r=0}^{L}$ partitions $S^1$ into gaps shorter than $|I|$, then we also have that $\{\Ra^{r}(0)\}_{|r|\leq \lceil \frac{L}{2}\rceil}$ will also partition $S^1$ into gaps shorter than $|I|$. Hence, for some $r\in\mathbb{Z}$ with $|r|\leq\lceil\frac{L}{2}\rceil$, we must have that $r\alpha \in I$, meaning $\rho(I)\leq |r|\leq\lceil\frac{L}{2}\rceil $.
    
    Finally, if we assume $|I|\leq \|q_{m}\alpha\|$, we note that if $a\alpha$ is an endpoint of $I$, and $r\in H(I)$, then $\|(r-a)\alpha\|<|I|\leq \|q_{m}\alpha\|$, then by the best approximation property of the convergents of $\alpha$, $|r-a|\geq q_{m+1}$, so $|r|\geq q_{m+1}-N$.
\end{proof}

We now give a general lemma bounding the growth of the series related to our metric.

 \begin{lem}\label{Lem:SeriesBound}
     For any $\lambda,\rho\ge 1$, and $\delta\in(0,1)$, there exists a uniform constant $c=c(\delta)>0$ such that 
    \begin{equation*}
    \rho^{-1/\delta}\leq\sum_{k\in\mathbb{Z}}\bigl(\rho+|k|\lambda\bigr)^{-1/\delta} \leq c\rho^{-1/\delta}\cdot\max\Bigl(1,\frac{\rho}{\lambda}\Bigr).
    \end{equation*}
 \end{lem}
\begin{proof}
    To begin with, we note that the lower bound is immediate from the fact that each term in the series is non-negative and $\rho^{-1/\delta}$ is itself a term in the series. Then, for the upper bound, we write
    \begin{align*}
        \sum_{k\in\mathbb{Z}}(\rho+|k|\lambda)^{-1/\delta} &= \lambda^{-1/\delta}\sum_{k\in\mathbb{Z}}\left(\frac{\rho}{\lambda}+|k|  \right)^{-1/\delta}.
    \end{align*}
     Let $B_{j}:=\{k\in\mathbb{Z}: j\cdot\frac{\rho}{\lambda} \leq |k|<(j+1)\cdot\frac{\rho}{\lambda}\}$ define a partition of $\mathbb{Z}$, and note the following. Firstly, $|B_{j}|\leq 2\cdot\frac{\rho}{\lambda}+2$. We also have that for each $k\in B_{j}$, $\frac{\rho}{\lambda}+|k| \geq (j+1)\cdot\frac{\rho}{\lambda}$, meaning $\left(\frac{\rho}{\lambda}+|k|\right)^{-1/\delta} \leq \left( (j+1)\cdot\frac{\rho}{\lambda}\right)^{-1/\delta}$. Then, we can write

    \begin{align*}
        \sum_{k\in\mathbb{Z}}(\rho+|k|\lambda)^{-1/\delta} &\leq \lambda^{-1/\delta}\sum_{j=0}^{\infty}\left(2\cdot\frac{\rho}{\lambda}+2\right)\left((j+1)\cdot\frac{\rho}{\lambda} \right)^{-1/\delta} \\ &= \lambda^{-1/\delta}\left(2\cdot\frac{\rho}{\lambda}+2\right)\sum_{j=0}^{\infty}\left((j+1)\cdot\frac{\rho}{\lambda} \right)^{-1/\delta} 
        \\ &= \lambda^{-1/\delta}\left(2\cdot\frac{\rho}{\lambda}+2\right)\sum_{j=1}^{\infty}\left(j\cdot\frac{\rho}{\lambda} \right)^{-1/\delta} \\ &= \lambda^{-1/\delta}\left(2\cdot\frac{\rho}{\lambda}+2\right)\left( \frac{\rho}{\lambda}\right)^{-1/\delta}\sum_{j=1}^{\infty}j^{-1/\delta} \\ &\leq \lambda^{-1/\delta}\cdot 4\max \left(1,\frac{\rho}{\lambda} \right)\cdot \left(\frac{\rho}{\lambda}\right)^{-1/\delta} \zeta(1/\delta) \\ &= 4\zeta(1/\delta) \cdot\rho^{-1/\delta} \cdot \max\left( 1,\frac{\rho}{\lambda}\right)
    \end{align*}
    So, writing $c=4\zeta(1/\delta)$ gives us the result.
\end{proof}
Combining the two lemmas above bounds the $\dOm$-length of an arc in terms of its first hitting time and the spacing of its hitting times, leading us to this corollary.

\begin{cor}\label{Cor:HittingSum}
Let $I\subset S^{1}$ be an open arc with first hitting time $\rho=\rho(I)$, and suppose that, for some $\sigma\in\mathbb{N}$ and $K\in\mathbb{N}$, every set of $\sigma$ consecutive integers contains at most $K$ elements of $H(I)$. Then
\begin{equation*}
    l_{\rho}\leq\sum_{r\in H(I)}l_{r} \leq c'K\,l_{\rho}\cdot\max\left(1,\frac{\rho}{\sigma}\right),
\end{equation*}
where $c'>0$ depends only on $\delta$.
\end{cor}
 
\begin{proof}
The lower bound holds since $\rho$ or $-\rho$ lies in $H(I)$, and $l_{-\rho}=l_{\rho}$. For the upper bound, note that every element of $H(I)$ has modulus at least $\rho$. For each $j\geq0$, at most $K$ positive elements of $H(I)$ lie in $[\rho+j\sigma,\rho+(j+1)\sigma)$, and each of these $r$ has $l_{r}\leq l_{\rho+j\sigma}$; the same holds for the absolute values of the negative elements. Hence, by Lemma \ref{Lem:SeriesBound},
\begin{align*}
    \sum_{r\in H(I)}l_{r}&\leq 2Kc_{\delta}\sum_{j\geq0}\big(\rho+j\sigma\big)^{-1/\delta}\\ &\leq 2Kc_{\delta}\sum_{j\in\mathbb{Z}}\big(\rho+|j|\sigma\big)^{-1/\delta}\\ &\leq2Kc_{\delta}c\rho^{-1/\delta}\cdot\max\left(1,\frac{\rho}{\sigma}\right).
\end{align*}
Since $1+\rho\leq2\rho$, we have $c_{\delta}\rho^{-1/\delta}\leq2^{1/\delta}l_{\rho}$, and the result follows with $c'=2^{1+1/\delta}c$.
\end{proof}
 
By Kac's Lemma (or equidistribution), consecutive hitting times of an arc $I$ are on average $\frac{1}{|I|}$ apart. The quantity $\sigma$ in Corollary \ref{Cor:HittingSum} controls the spacing of the return times rather than their mean, and the factor $\max\left(1,\frac{\rho}{\sigma} \right) $ measures how late the first hit is compared to this spacing. We first treat the arcs for which this factor is bounded.

\begin{lem}\label{lem:technicalLem}
    Let $(0,n\alpha)\subset S^{1}$ for $n\in\mathbb{Z}$ be a subinterval of the circle, let $n_{0}\in\mathbb{N}$ be the smallest natural number such that $q_{n_{0}}\alpha\in [0,n\alpha]$, and let $\rho =\rho((0,n\alpha))\in\mathbb{N}$ be its first hitting time. Suppose that either
    \begin{equation*}
    \|n\alpha\|=\|q_{n_{0}}\alpha\| \;\text{ or }\; \|n\alpha\|=\|(q_{n_{0}}-q_{n_{0}-1})\alpha\|.
    \end{equation*}
    Then  
    \begin{equation*}
        \dOm(h^{-1}(0),h^{-1}(n\alpha)) = \Theta(l_{\rho}),
    \end{equation*}
    where the implied constants depend only on $\delta$.
\end{lem}
\begin{proof}
    When considering preimages under our semiconjugacy $h$, since the points $0,n\alpha\in \{k\alpha\bmod 1:k\in\mathbb{Z}\}$ by definition, the definition of $h$ tells us that their preimages will be the doubleton sets $\{ g_{r}(0),g_{l}(0)\}$ and $\{g_{r}(n\alpha), g_{l}(n\alpha)\}$ respectively. To make our arguments formal, we will choose the preimages such that the interval formed by these preimages has minimal length. In other words, if we have the interval $(0,n\alpha)\subset S^{1}$ oriented clockwise, we take $(h^{-1}(0),h^{-1}(n\alpha))=(g_{r}(0),g_{l}(n\alpha))$. We note that taking preimages in this way naturally leads to our $\dOm$ metric measuring our intervals by the hitting times of open intervals, in particular, excluding the endpoints as potential hitting times. Writing $I=(0,n\alpha)$, we have that
    \begin{equation*}
         \dOm(h^{-1}(0),h^{-1}(n\alpha)) =\sum_{r\in H(I)} l_{r}.
    \end{equation*}
    By Corollary \ref{Cor:HittingSum} with $K=1$, it suffices to show that any two distinct elements of $H(I)$ differ by at least $\rho$. Before we start, we recall the following facts. First, that the points $q_{k}\alpha$ fall on either side of $0$, then that $\|q_{k-1}\alpha\|=a_{k+1}\|q_{k}\alpha\|+\|q_{k+1}\alpha\|$, and finally that by the best approximation property of convergents, $\|r\alpha\|\geq \|q_{k-1}\alpha\|$ whenever $0<|r|<q_{k}$. We treat the two cases from the statement individually.

    $\underline{\|n\alpha\|=\|q_{n_{0}}\alpha\|}$: In this case, $|I|=\|q_{n_{0}}\alpha\|$, so by Lemma \ref{Lem:HittingSeparation}$(i)$, setting $m=n_{0}+1$, any two distinct elements of $H(I)$ differ by at least $q_{n_{0}+1}$. We now show that $\rho=q_{n_{0}+1}$. The point $-q_{n_{0}+1}\alpha$ lies on the same side of $0$ as $q_{n_{0}}\alpha$, at distance $\|q_{n_{0}+1}\alpha\|<\|q_{n_{0}}\alpha\|$, so $-q_{n_{0}+1}\in H(I)$. On the other hand, if $0<|r|<q_{n_{0}+1}$, then $\|r\alpha\|\geq \|q_{n_{0}}\alpha\|=|I|$, so $r\notin H(I)$. Hence $\rho=q_{n_{0}+1}$.

    $\underline{\|n\alpha\|=\|(q_{n_{0}}-q_{n_{0}-1})\alpha\|}$: Since $q_{n_{0}}\alpha,q_{n_{0}-1}\alpha$ fall on opposite sides of $0$, in this case, $|I|=\|q_{n_{0}}\alpha\|+\|q_{n_{0}-1}\alpha\|$, and $I$ lies on the same side of $0$ as $q_{n_{0}}\alpha$. The point $-q_{n_{0}-1}\alpha$ lies on the same side of $0$, at distance $\|q_{n_{0}-1}\alpha\|<|I|$, and so $-q_{n_{0}-1}\in H(I)$. We note that if $n_{0}=1$, then $q_{n_{0}-1}=1$, so $\rho=1$, and the distinct hitting times trivially differ by at least 1, so we may therefore assume that $n_{0}\geq 2$. We have that, if $0<|r|<q_{n_{0}-1}$, then $\|r\alpha\|\geq \|q_{n_{0}-2}\alpha\|=a_{n_{0}}\|q_{n_{0}-1}\alpha\|+\|q_{n_{0}}\alpha\|\geq |I|$, so $r\notin H(I)$. Hence $\rho=q_{n_{0}-1}$. Then, since $|I|\leq \|q_{n_{0}-2}\alpha\|$, by Lemma \ref{Lem:HittingSeparation}$(i)$ with $m=n_{0}-1$, we have that two distinct elements of $H(I)$ must differ by at least $q_{n_{0}-1}=\rho$.

    In each case, we have shown that the hitting times of $I$ differ by at least $\rho$, so Corollary \ref{Cor:HittingSum} gives $l_{\rho}\leq  \dOm(h^{-1}(0),h^{-1}(n\alpha))  \leq c'l_{\rho}$.
\end{proof}

We now strengthen this result to apply to intervals around the origin and across $S^{1}$. We are still interested in intervals formed by orbits of $\Ra(0)$.

\begin{cor}\label{cor:BigAppendixCor}
   Let $\alpha\in(0,1)$ be irrational. Let $(n\alpha,m\alpha)\subset S^{1}$ for $n,m\in\mathbb{Z}$ be a subinterval of the circle such that 
    \begin{equation*}
        \|(m-n)\alpha\|=\|q_{n_{0}}\alpha\| \; \text{ or }\; \|(m-n)\alpha\|=\|(q_{n_{0}}-q_{n_{0}-1})\alpha\|,
    \end{equation*}
    for some $n_{0}\in\mathbb{N}$, with $n_{0}\geq 2$ in the second case. Let $\rho=\rho((n\alpha,m\alpha))$ be the first hitting time by $(n\alpha,m\alpha)$ of the orbit of $\Ra(0)$.
     Then,
    \begin{equation*}
        \dOm(h^{-1}(n\alpha),h^{-1}(m\alpha)) = \Theta(l_{\rho}),
    \end{equation*}
    where the implied constants depend only on $\delta$.
\end{cor}
\begin{proof}
    Write $I=(n\alpha,m\alpha)$. Choosing preimages as in the proof of Lemma \ref{lem:technicalLem}, we have $ \dOm(h^{-1}(n\alpha),h^{-1}(m\alpha)) =\sum_{r\in H(I)}l_{r}$. By Corollary \ref{Cor:HittingSum}, it suffices to find $\sigma, K$ as in the statement with $K\leq 2$ and $\rho\leq 2\sigma$.

    $\underline{\|(m-n)\alpha\|=\|q_{n_{0}}\alpha\|}$: In this case,  $|I|=\|q_{n_{0}}\alpha\|$, so by Lemma \ref{Lem:HittingSeparation}$(i)$, applied with $n_{0}+1$ in place of $m$, we can take $\sigma=q_{n_{0}+1}$ and $K=1$. By Theorem \ref{Thm:vanRavenstein}, the points $\{\Ra^{r}(0)\}_{r=0}^{q_{n_{0}}+2q_{n_{0}+1}}$ partition $S^1$ into gaps of length less than $\|q_{n_{0}}\alpha\|$. If we take $L=q_{n_{0}}+2q_{n_{0}+1}$ and $M=\lceil\frac{1}{2}(q_{n_{0}}+2q_{n_{0}+1})\rceil$, Lemma \ref{Lem:BoundsOnRho} gives us that $\rho\leq M\leq 2q_{n_{0}+1}$.

    $\underline{\|(m-n)\alpha\|=\|(q_{n_{0}}-q_{n_{0}-1})\alpha\|}$: In this case, $|I|=\|q_{n_{0}}\alpha\|+\|q_{n_{0}-1}\alpha\|$, and so by Lemma \ref{Lem:HittingSeparation}$(ii)$ taking $m=n_{0}$, we may take $\sigma=q_{n_{0}}$ and $K=2$. By Theorem \ref{Thm:vanRavenstein}, the points in $\{\Ra^{r}(0)\}_{r=0}^{q_{n_{0}}+q_{n_{0}-1}}$ partition $S^1$ into gaps of lengths $\|q_{n_{0}-1}\alpha\|$ and $\|q_{n_{0}}\alpha\|$, both of which are shorter than $|I|$. Taking $L=q_{n_{0}}+q_{n_{0}-1}$ and $M=\lceil \frac{1}{2}(q_{n_{0}}+q_{n_{0}-1})\rceil$ gives $\rho\leq M\leq q_{n_{0}}$, again by Lemma \ref{Lem:BoundsOnRho}.

    In both cases, $\max\left(1,\frac{\rho}{\sigma}\right)\leq 2$, and $K\leq 2$, so Corollary \ref{Cor:HittingSum} gives us that
    \begin{equation*}
        l_{\rho}\leq \dOm(h^{-1}(n\alpha),h^{-1}(m\alpha)) \leq 4c'l_{\rho}.
    \end{equation*}
\end{proof}

\begin{rmk}
    For the two interval types covered above, the first hitting time is at most a bounded multiple of the spacing $\sigma$ between consecutive hitting times, so the factor $\max\left( 1,\frac{\rho}{\sigma}\right)$ in Corollary \ref{Cor:HittingSum} is bounded, and the $\dOm$--length is $\Theta(l_{\rho})$. This is not the case for the complementary intervals of $q_{n_{0}+1,i-1}$-type, with $i\geq 2$, considered in Proposition \ref{Prop:BigAppendixProp}. These have length less than $\|q_{n_{0}-1}\alpha\|$, so by Lemma \ref{Lem:HittingSeparation} (i), their hitting times are spaced at least $q_{n_{0}}$ apart, while their first hitting times are $\Theta(iq_{n_{0}})$, meaning $\frac{\rho}{\sigma}=\Theta(i)$. These intervals are dealt with in Proposition \ref{Prop:BigAppendixProp}.
\end{rmk}

\subsection{Lengths of complementary intervals}
Before proving the upcoming proposition, we will need a technical lemma.

\begin{lem}\label{Lem:FinalAppendixLem}
    Let $\delta\in (0,1)$ and let $2\leq i \leq a_{n_{0}+1}-1$. Then
    \begin{equation*}
        \sum_{j=i+1}^{a_{n_{0}+1}}j^{-1/\delta} = \Theta \left( i^{1-1/\delta}\cdot \frac{a_{n_{0}+1}-(i-1)}{a_{n_{0}+1}}\right)
    \end{equation*}
    with implied constants depending only on $\delta$, uniformly in $i$ and $a_{n_{0}+1}$.
\end{lem}
\begin{proof}
    We first note that since $a_{n_{0}+1}-i\geq 1$, the sum is non-empty, and we have that
    \begin{equation}\label{eqn:appendixFinalLemAn0,i inequalities}
        a_{n_{0}+1}-i\leq a_{n_{0}+1}-(i-1)\leq 2(a_{n_{0}+1}-i)
    \end{equation}
    with the second inequality being a consequence of the fact that $a_{n_{0}+1}-i\geq 1$. We separate this argument into two cases.
    
    \underline{Case 1: $i\geq\frac{a_{n_{0}+1}}{2}$.} In this case, every $j$ appearing in the sum satisfies that $i<j\leq a_{n_{0}+1}\leq 2i$, so we get that
    \begin{equation*}
        (2i)^{-1/\delta} \leq j^{-1/\delta} \leq i^{-1/\delta}.
    \end{equation*}
    Applying this to the $a_{n_{0}+1}-i$ terms in the series, we get that
    \begin{equation*}
        2^{-1/\delta}(a_{n_{0}+1}-i)i^{-1/\delta} \leq \sum_{j=i+1}^{a_{n_{0}+1}}j^{-1/\delta}\leq (a_{n_{0}+1}-i)i^{-1/\delta}.
    \end{equation*}
    Then, by \eqref{eqn:appendixFinalLemAn0,i inequalities}, we deduce that $\sum_{j=i+1}^{a_{n_{0}+1}}j^{-1/\delta}=\Theta((a_{n_{0}+1}-(i-1))i^{-1/\delta})$. To see that this does in fact match the statement, we note that since $\frac{i}{2}<\frac{a_{n_{0}+1}}{2}\leq i$, we have $\frac{1}{2}\leq \frac{i}{a_{n_{0}+1}}<1$. This means that
    \begin{align*}
        i^{1-1/\delta}\cdot\frac{a_{n_{0}+1}-(i-1)}{a_{n_{0}+1}} &= \frac{i}{a_{n_{0}+1}} \cdot (a_{n_{0}+1}-(i-1))\cdot i^{-1/\delta} \\ &= \Theta\left((a_{n_{0}+1}-(i-1))i^{-1/\delta} \right).
    \end{align*}
    So, the claim holds when $i\geq \frac{a_{n_{0}+1}}{2}$.

    \underline{Case 2: $i<\frac{a_{n_{0}+1}}{2}$.} We begin by noting that when $i<\frac{a_{n_{0}+1}}{2}$, we have that
    \begin{equation*}
        \frac{1}{2}<\frac{a_{n_{0}+1}-(i-1)}{a_{n_{0}+1}}\leq 1.
    \end{equation*}
    This means that it is sufficient to show that the sum is of order $\Theta\left(i^{1-1/\delta}\right)$. For the lower bound, we note that since $i<\frac{a_{n_{0}+1}}{2}$, we get that $2i<a_{n_{0}+1}$, which means that each of $j=i+1,i+2,\dots,2i$ appear in the sum. This means that
    \begin{align*}
        \sum_{j=i+1}^{a_{n_{0}+1}}j^{-1/\delta} &\geq \sum_{j=i+1}^{2i}j^{-1/\delta} \\ &\geq i\cdot (2i)^{-1/\delta} \\ &= 2^{-1/\delta} \cdot i^{1-1/\delta}.
    \end{align*}
    For the upper bound, we extend the sum for all $j>i$, and substitute $j=i+k$. Then
    \begin{align*}
        \sum_{j=i+1}^{a_{n_{0}+1}}j^{-1/\delta} &\leq \sum_{k\geq 1}(i+k)^{-1/\delta} \\ &\leq \sum_{k\in\mathbb{Z}}(i+|k|)^{-1/\delta}.
    \end{align*}
    Applying Lemma \ref{Lem:SeriesBound} setting $\rho= i$ and $\lambda=1$, we have that $\max\left(1,\frac{\rho}{\lambda}\right)=i$, and so there exists some $c=c(\delta)>0$ such that
    \begin{equation*}
        \sum_{k\in\mathbb{Z}}(i+|k|)^{-1/\delta}\leq c\cdot i^{1-1/\delta},
    \end{equation*}
    giving us the result.
\end{proof}

\begin{prop}\label{Prop:BigAppendixProp}
    Let $\alpha\in (0,1)$ be irrational and $\delta\in(0,1)$. Let $n_{0}\ge 2$ and $2\leq i \leq a_{n_{0}+1}$, and let $N\in \mathbb{N}$ satisfy $q_{n_{0}+1,i-1}< 2N+1 \leq q_{n_{0}+1,i}$. Let $I\subset S^1$ be some complementary interval of $\mathcal{J}_{N}$ whose first hitting time is $\rho=\rho(I)$. Then, 
    \begin{equation*}
        |I|_{\dOm}  = \Theta(l_{\rho}\cdot\max (1,\rho|I|))
    \end{equation*}
    where $|I|$ denotes the Lebesgue-length of the interval $I$ in $S^1$, while $|I|_{\dOm}$ is the $\dOm$-length of $I$.
    In particular, when $I$ is a $q_{n_{0}}$-type interval, then $\max(1,\rho|I|)=\Theta(1)$ and
        \begin{equation*}
            |I|_{\dOm} = \Theta(l_{q_{n_{0}+1}})=\Theta(q_{n_{0}+1}^{-1/\delta}).
        \end{equation*}
    If $I$ is a $q_{n_{0}+1,i-1}$-type interval, then $\rho=\Theta(iq_{n_{0}})$, $|I|=\Theta\left((a_{n_{0}+1}-(i-1))q_{n_{0}+1}^{-1} \right)$ and 
    \begin{equation*}
        \max(1,\rho|I|)=\Theta\left(\frac{i(a_{n_{0}+1}-(i-1))}{a_{n_{0}+1}} \right).
        \end{equation*}
        In each case, the implied constants depend only on $\delta$.
\end{prop}
\begin{proof}
    We begin by listing the known information about the interval types present in $\mathcal{J}_{N}$ with $N$ as described in the statement of this proposition. By Corollary \ref{Cor:Integer3GapThrm}, the complementary intervals of $\J_{N}$ are either $q_{n_{0}}$-type, $q_{n_{0}+1,i-1}$-type, or $(q_{n_{0}}-q_{n_{0}+1,i-1})$-type. We note that the last of these types coincides with the $q_{n_{0}+1,i-2}$-type intervals: indeed, $\|(q_{n_{0}}-q_{n_{0}+1,i-1})\alpha\|=\|q_{n_{0}+1,i-2}\alpha\|$ since $q_{n_{0}+1,i-1}=(i-1)q_{n_{0}}+q_{n_{0}-1}$ gives $q_{n_{0}+1,i-2}=q_{n_{0}+1,i-1}-q_{n_{0}}$. The argument below for the $q_{n_{0}+1,i-1}$-type intervals applies to the $q_{n_{0}+1,i-2}$-type intervals with $i-1$ in place of $i$ (when $i=2$, these are $q_{n_{0}-1}$-type intervals, covered by Corollary \ref{cor:BigAppendixCor}), so it suffices to treat the first two interval types.

    \underline{$I\in\mathcal{I}_{q_{n_{0}}}$:} Since $2N+1\leq q_{n_{0}+1}$, Lemma \ref{Lem:BoundsOnRho}$(iii)$ with $m=n_{0}$ gives us that $\rho>\frac{q_{n_{0}+1}}{2}$, and Lemma \ref{Lem:BoundsOnRho}$(ii)$ with $L=q_{n_{0}}+2q_{n_{0}+1}$ as in the proof of Corollary \ref{cor:BigAppendixCor}, gives us that $\rho\leq 2q_{n_{0}+1}$. So, $\rho =\Theta(q_{n_{0}+1})$, and $|I|=\|q_{n_{0}}\alpha\|=\Theta(q_{n_{0}+1}^{-1})$, so 
    \begin{equation*}
        \rho |I|= \Theta(1).
    \end{equation*}
    Corollary \ref{cor:BigAppendixCor} applies to $I$, since any $q_{n_{0}}$-type interval as described in the statement of this proposition must indeed have endpoints $(n\alpha,m\alpha)$ for some $n,m\in\mathbb{Z}$ with $\|(m-n)\alpha\|=\|q_{n_{0}}\alpha\|$. Hence, we have that
    \begin{align*}
        |I|_{\dOm}=\Theta(l_{\rho}\cdot\max(1,\rho|I|)),
    \end{align*}
    as required.

    \underline{$I\in\mathcal{I}_{q_{n_{0}+1,i-1}}$:} We begin by noting that 
    \begin{equation*}
        \|q_{n_{0}-1}\alpha\|=a_{n_{0}+1}\|q_{n_{0}}\alpha\| +\|q_{n_{0}+1}\alpha\|.
    \end{equation*}
     Since $q_{n_{0}+1,i-1}=q_{n_{0}-1}+(i-1)q_{n_{0}}$ and the errors of $q_{n_{0}}\alpha$ and $q_{n_{0}-1}\alpha$ have opposite signs, each increase of index by $q_{n_{0}}$ decreases the norm by $\|q_{n_{0}}\alpha\|$ so that
    \begin{equation*}
        \|q_{n_{0}+1,i-1}\alpha\|= (a_{n_{0}+1}-(i-1))\|q_{n_{0}}\alpha\| +\|q_{n_{0}+1}\alpha\|.
    \end{equation*}
      Since $\|q_{n_{0}+1}\alpha\|<\|q_{n_{0}}\alpha\|$ and, by \eqref{eqn:CFnorminequality2}, $\|q_{n_{0}}\alpha\|=\Theta(q_{n_{0}+1}^{-1})$, we have that 
    \begin{equation*}
     |I|=\Theta\left(\frac{a_{n_{0}+1}-(i-1)}{q_{n_{0}+1}}\right) \text{ for } I\in \mathcal{I}_{q_{n_{0}+1,i-1}}.
    \end{equation*}
     Since $I$ is a complementary interval of $\J_{N}$, we must have that $\rho>N\geq\frac{1}{2}q_{n_{0}+1,i-1}$, since $2N+1>q_{n_{0}+1,i-1}$. By Theorem \ref{Thm:vanRavenstein}, the points $r\alpha$ with $0\leq r< q_{n_{0}+1,i}+q_{n_{0}}$ partition $S^1$ into gaps of length $\|q_{n_{0}}\alpha\|$ and $\|q_{n_{0}+1,i}\alpha\|$, both of which are shorter than $|I|$, so Lemma \ref{Lem:BoundsOnRho}$(ii)$ tells us that $\rho\leq \lceil\frac{1}{2}(q_{n_{0}+1,i}+q_{n_{0}})\rceil \leq q_{n_{0}+1,i} $. Since $jq_{n_{0}}\leq q_{n_{0}+1,j}\leq (j+1)q_{n_{0}} $ for $j\geq 1$, we obtain
    \begin{equation*}
        \rho=\Theta(iq_{n_{0}}).
    \end{equation*}
    From this, first noting that $q_{n_{0}+1}=\Theta(a_{n_{0}+1}q_{n_{0}})$, we deduce that 
    \begin{equation*}
        \rho|I|=\Theta\left(\frac{i(a_{n_{0}+1}-(i-1))}{a_{n_{0}+1}}\right).
    \end{equation*}
    When $i=a_{n_{0}+1}$, the product $i(a_{n_{0}+1}-(i-1))$ equals $a_{n_{0}+1}$; and when $2\leq i\leq a_{n_{0}+1}-1$, both terms are at least $2$, so the product is at least $a_{n_{0}+1}+1$. In either case, $i(a_{n_{0}+1}-(i-1))\geq a_{n_{0}+1}$, and so $\rho|I|$ is bounded below by a positive constant. Hence $\max(1,\rho|I|)=\Theta(\rho|I|)$ for the remainder of this case, and it will suffice to prove that 
    \begin{align*}
        |I|_{\dOm}=\Theta\left(l_{\rho}\cdot\frac{i(a_{n_{0}+1}-(i-1))}{a_{n_{0}+1}} \right).
    \end{align*}
    We will first deal with the edge case where $i=a_{n_{0}+1}$. Here, we have that $|I|=\|q_{n_{0}}\alpha\|+\|q_{n_{0}+1}\alpha\|=\|(q_{n_{0}+1}-q_{n_{0}})\alpha\|$, so by Corollary \ref{cor:BigAppendixCor}, applied with $n_{0}$ replaced by $n_{0}+1$, we obtain $|I|_{\dOm}=\Theta(l_{\rho})$. At $i=a_{n_{0}+1}$ the quantity $\frac{i(a_{n_{0}+1}-(i-1))}{a_{n_{0}+1}}$ equals $1$, so $\max(1,\rho|I|)=\Theta(1)$ and the result holds in this case. So we now restrict to $2\leq i\leq a_{n_{0}+1}-1$.
    
    We now want to understand the hitting set, which we denote by $H=H(I)=\{r\in\mathbb{Z}: r\alpha\in I\}$, in order to set up the series defining the $\dOm$-length of $I$. By Corollary \ref{Cor:Integer3GapThrm}, for each $n_{0}$, there are at most $q_{n_{0}}$ $q_{n_{0}+1,i-1}$-type intervals, and each is split exactly once, forming another $q_{n_{0}}$-type interval and a $q_{n_{0}+1,i}$-type interval. This cycle will repeat, creating $q_{n_{0}}$ more $q_{n_{0}}$-type intervals and a $q_{n_{0}+1,i+1}$-type, and so on, until we have partitioned $I$ into $a_{n_{0}+1}-i$ $q_{n_{0}}$-type intervals and a $(q_{n_{0}+1}-q_{n_{0}})$-type interval. For $1\leq k\leq a_{n_{0}+1}-i$, let $r_{k}$ denote the hitting time at which $I$ is split for the $k^{th}$ time. The piece of $I$ remaining after $k-1$ splits is a $q_{n_{0}+1,i+k-2}$-type complementary interval of $\J_{|r_{k-1}|}$, where $r_{0}=N$, and $r_{k}$ is its first hitting time. By Corollary \ref{Cor:Integer3GapThrm}, this type only occurs in $\J_{N'}$ when $2N'+1>q_{n_{0}+1,i+k-2}$, so the argument above with $i+k-1$ in place of $i$ gives us that $|r_{k}|\leq \lceil  \frac{1}{2}q_{n_{0}+1,i+k}\rceil \leq \lceil\frac{1}{2}q_{n_{0}+1}\rceil$, and 
    \begin{equation}\label{eqn:|rk|=Theta((i+k)qn0}
        |r_{k}|=\Theta((i+k)q_{n_{0}}).
    \end{equation}
    Since we have already dealt with the $q_{n_{0}}$-type intervals and the $(q_{n_{0}+1}-q_{n_{0}})$-type intervals in Corollary \ref{cor:BigAppendixCor}, we can use this to separate our series describing the $\dOm$-length of $I$ into the hitting times $r_{k}$ and the subsequent $\dOm$-lengths of the remaining pieces. Let $\mathcal{P}$ denote the collection of complementary intervals into which $I$ is partitioned by this process, consisting of the $a_{n_{0}+1}-i$ intervals of $q_{n_{0}}$-type, along with the single terminal interval, which is of $(q_{n_{0}+1}-q_{n_{0}})$-type, giving us that $|\mathcal{P}|=a_{n_{0}+1}-(i-1)$. 

    Since $\sum_{n\in\mathbb{Z}}l_{n}=1$, the $\dOm$-length of an interval is precisely the sum of the lengths $l_{r}$ of the Denjoy intervals it contains. Now, $I$ is the disjoint union of the intervals of $\mathcal{P}$, together with the Denjoy intervals $J_{r_{1}},\dots , J_{r_{a_{n_{0}+1}-i}}$ inserted at the splitting times, so every $r$ with $r\alpha\in I$ lies in exactly one member of $\mathcal{P}$ or is one of the $r_{k}$. In other words, we can write
    \begin{equation*}
        |I|_{\dOm}= \sum_{k=1}^{a_{n_{0}+1}-i} l_{r_{k}}+\sum_{J\in \mathcal{P}} |J|_{\dOm} .
    \end{equation*}

    To begin with, we look at the hitting times represented in the first of these sums. Using \eqref{eqn:|rk|=Theta((i+k)qn0}, we have that 
    \begin{align*}
        \sum_{k=1}^{a_{n_{0}+1}-i} l_{r_{k}} = \Theta\left(\sum_{k=1}^{a_{n_{0}+1}-i}((i+k)q_{n_{0}})^{-1/\delta}  \right).
    \end{align*}
    Changing variables to write $j=i+k$, this becomes
    \begin{align*}
        \sum_{k=1}^{a_{n_{0}+1}-i} l_{r_{k}} &= \Theta\left(\sum_{j=i+1}^{a_{n_{0}+1}}(jq_{n_{0}})^{-1/\delta}  \right) \\&= \Theta\left(q_{n_{0}}^{-1/\delta}\sum_{j=i+1}^{a_{n_{0}+1}} j^{-1/\delta} \right).
    \end{align*}
    By Lemma \ref{Lem:FinalAppendixLem}, we have that
    \begin{align*}
        \sum_{j=i+1}^{a_{n_{0}+1}}j^{-1/\delta} =\Theta\left(i^{1-1/\delta}\cdot\frac{a_{n_{0}+1}-(i-1)}{a_{n_{0}+1}} \right),
    \end{align*}
    so combining this with the above gives us that the first sum is
    \begin{align*}
        \Theta \left(i\cdot(iq_{n_{0}})^{-1/\delta}\cdot\frac{a_{n_{0}+1}-(i-1)}{a_{n_{0}+1}} \right) = \Theta\left(l_{\rho}\cdot\frac{i(a_{n_{0}+1}-(i-1))}{a_{n_{0}+1}} \right).
    \end{align*}
    For the remaining terms in the series, we note that the members of $\mathcal{P}$ are of $q_{n_{0}}$-type or, in the case of the terminal interval, of $(q_{n_{0}+1}-q_{n_{0}})$-type. For the former,  each such $J$ is a complementary interval of $\J_{|r_{k}|}$ for some $k$, so Lemma \ref{Lem:BoundsOnRho}$(iii)$ and the bound $|r_{k}|\leq \lceil\frac{1}{2}q_{n_{0}+1}\rceil $ give $\rho(J)\geq q_{n_{0}+1}-\lceil\frac{1}{2}q_{n_{0}+1}\rceil$, while Lemma \ref{Lem:BoundsOnRho}$(ii)$ gives $\rho(J)\leq 2q_{n_{0}+1}$. Hence $\rho(J)=\Theta(q_{n_{0}+1})$, and Corollary \ref{cor:BigAppendixCor} applies, giving $|J|_{\dOm}=\Theta\left(q_{n_{0}+1}^{-1/\delta}\right)$.
    
    For the latter, its length is $\|(q_{n_{0}+1}-q_{n_{0}})\alpha\|=\|q_{n_{0}}\alpha\|+\|q_{n_{0}+1}\alpha\|$, and its first hitting time is $\Theta(q_{n_{0}+1})$ by the argument giving us \eqref{eqn:|rk|=Theta((i+k)qn0} with $k=a_{n_{0}+1}-i+1$, so Corollary \ref{cor:BigAppendixCor} applies with $n_{0}$ replaced with $n_{0}+1$, as in the edge case $i=a_{n_{0}+1}$ treated above, giving $O\left(q_{n_{0}+1}^{-1/\delta} \right)$. Since $|\mathcal{P}|=a_{n_{0}+1}-(i-1)$, the remainder of the sum is $O\left((a_{n_{0}+1}-(i-1))q_{n_{0}+1}^{-1/\delta} \right)$, and since $a_{n_{0}+1}-i\geq 1$ for $i<a_{n_{0}+1}$, the $q_{n_{0}}$-type members alone give a matching lower bound of the same order. Hence, the remainder of the sum is of order
    \begin{equation*}
        \Theta\left((a_{n_{0}+1}-(i-1))q_{n_{0}+1}^{-1/\delta} \right).
    \end{equation*}
    While this does not appear to be in the form in the statement, by considering their ratio, we see that 
    \begin{align*}
        \frac{(a_{n_{0}+1}-(i-1))q_{n_{0}+1}^{-1/\delta}}{\left(l_{\rho}\cdot\frac{i(a_{n_{0}+1}-(i-1))}{a_{n_{0}+1}}\right)} &= \frac{a_{n_{0}+1}}{i\cdot l_{\rho}\cdot q_{n_{0}+1}^{1/\delta}} \\ &= \Theta \left(\frac{a_{n_{0}+1}}{i} \cdot \left(\frac{iq_{n_{0}}}{q_{n_{0}+1}} \right)^{1/\delta}\right).
    \end{align*}
    Then, since $q_{n_{0}+1}=\Theta(a_{n_{0}+1}q_{n_{0}})$, we deduce that 
    \begin{align*}
        \Theta\left( \frac{(a_{n_{0}+1}-(i-1))q_{n_{0}+1}^{-1/\delta}}{\left(l_{\rho}\cdot\frac{i(a_{n_{0}+1}-(i-1))}{a_{n_{0}+1}}\right)}\right) &= \Theta\left(\left( \frac{i}{a_{n_{0}+1}}\right)^{1/\delta-1} \right)\\&= O(1) 
    \end{align*}
    since both $i\leq a_{n_{0}+1}$ and $1/\delta >1$. Thus, we have shown that 
    \begin{equation*}
        |I|_{\dOm} = \Theta\left(l_{\rho}\cdot \frac{i(a_{n_{0}+1}-(i-1))}{a_{n_{0}+1}}\right).
    \end{equation*}    
\end{proof}

\subsection{Intervals of the same type}
This subsection aims to prove a lemma stating that intervals of the same type are 'essentially' of equal length. 

\begin{lem}\label{BigLem}
    Let $\alpha\in (0,1)\setminus\mathbb{Q}$. Let $f$ be a classical Denjoy map with rotation number $\alpha$ and minimal set $\OmegaAD$. Then, for each $ N\in\mathbb{N}$ with $2N+1>q_{2}$, $\mathcal{J}_{N}$ is made up of at most 3 different interval types, each of which has the property that intervals of a given type are of equal $\dOm$-length up to some uniformly bounded multiplicative constant. In other words, there exists $C=C(\delta)>0$ such that, for every $N$, every interval type in $\J_{N}$, and every pair of intervals $I,J$ of that type, 
    \begin{equation*}
        C^{-1}\leq \frac{|I|_{\dOm}}{|J|_{\dOm}}\leq C.
    \end{equation*}
\end{lem}
\begin{proof}
    The idea of the proof is that two intervals of the same type, say $x$-type,  have the same length in $S^1$, namely $\|x\alpha\|$. By Corollary \ref{cor:BigAppendixCor} and Proposition \ref{Prop:BigAppendixProp}, their $\dOm$--lengths are the same function of their first hitting time $\rho$, up to constants depending only on $\delta$. That function changes by a factor of at most $M^{1+1/\delta}$ when $\rho$ changes by a factor of $M$, so this lemma will come down to showing that for each interval type, every first hitting time lies in $(T,MT]$ for a single $T$. In particular, we will show that taking $M=4$ is sufficient.

    Fix $N$ with $2N+1>q_{2}$, and let $n_{0}\geq 2$ satisfy $q_{n_{0}}<2N+1\leq q_{n_{0}+1}$, so that $q_{n_{0}}\leq 2N<q_{n_{0}+1}$. As in Corollary \ref{Cor:Integer3GapThrm}, we assume that $n_{0}+1$ is odd, the even case being symmetric. By Corollary \ref{Cor:Integer3GapThrm}, $\J_{N}$ has at most three interval types, and as in the statement of the Corollary, we separate the two cases
    \begin{enumerate}[(i)]
        \item $q_{n_{0}}<2N+1\leq q_{n_{0}+1,1}$  
        \item $q_{n_{0}+1,i-1}<2N+1\leq q_{n_{0}+1,i}$ for some $2\leq i\leq a_{n_{0}+1}$.
    \end{enumerate}

    Let $I,J$ be two intervals of the same type, say $x$-type, so that $|I|=|J|=\|x\alpha\|$. In case $(i)$, the three possible interval types are $q_{n_{0}}$-type, $q_{n_{0}-1}$-type, and $(q_{n_{0}}-q_{n_{0}-1})$-type, all of which satisfy $|I|_{\dOm}=\Theta(l_{\rho(I)})$ by Corollary \ref{cor:BigAppendixCor}. In case $(ii)$, Proposition \ref{Prop:BigAppendixProp} applies, giving $|I|_{\dOm}=\Theta\left(l_{\rho(I)}\max(1,\rho(I)\|x\alpha\|) \right)$. The same holds for $J$. We note that if $\rho(I)<4\rho(J)$ and $\rho(J)<4\rho(I)$, then $l_{\rho(I)}\leq 4^{1/\delta}l_{\rho(J)}$ and $\max\{1,\rho(I)\|x\alpha\|\}\leq 4\max\{1,\rho(J)\|x\alpha\|\}$, so in either case,  
    \begin{equation*}
        |I|_{\dOm}\leq c\cdot 4^{1+1/\delta}|J|_{\dOm}
    \end{equation*}
    where $c$ depends only on $\delta$. Hence, by symmetry, it suffices to show that, for each interval type $x$ in $\J_{N}$, there is some $T>0$ such that 
    \begin{equation}\label{Eqn:T<rho(I)<4T}
        T<\rho(I)\leq 4T \;\; \text{ for every }x\text{-type interval } I,
    \end{equation}
    with the lemma following taking $C=c4^{1+1/\delta}$.

    We now examine each interval type using Lemma \ref{Lem:BoundsOnRho}. The $q_{n_{0}}$-type intervals have $|I|=\|q_{n_{0}}\alpha\|$, so, since $2N+1\leq q_{n_{0}+1}$, Lemma \ref{Lem:BoundsOnRho}$(iii)$ with $m=n_{0}$ gives $\rho(I)>\frac{q_{n_{0}+1}}{2}$, and Lemma \ref{Lem:BoundsOnRho}$(ii)$ with $L=q_{n_{0}}+2q_{n_{0}+1}$, as in the proof of Corollary \ref{cor:BigAppendixCor}, gives $\rho(I)\leq \lceil\frac{1}{2}(q_{n_{0}}+2q_{n_{0}+1})\rceil<2q_{n_{0}+1}$, so \eqref{Eqn:T<rho(I)<4T} is satisfied with $T=\frac{q_{n_{0}+1}}{2}$. 
    
    The remaining intervals in case $(i)$, of $q_{n_{0}-1}$-type and $(q_{n_{0}}-q_{n_{0}-1})$-type, have length at least $\|q_{n_{0}-1}\alpha\|$. As in the proof of Corollary \ref{cor:BigAppendixCor}, the points $r\alpha$ with $0\leq r<q_{n_{0}-1}+2q_{n_{0}}$ partition $S^1$ into gaps shorter than $\|q_{n_{0}-1}\alpha\|$, so Lemma \ref{Lem:BoundsOnRho}$(ii)$ gives $\rho(I)\leq \lceil\frac{1}{2}(q_{n_{0}-1}+2q_{n_{0}})\rceil<2q_{n_{0}}$, and since $\rho(I)>N\geq\frac{q_{n_{0}}}{2}$, by Lemma \ref{Lem:BoundsOnRho}$(i)$, \eqref{Eqn:T<rho(I)<4T} holds with $T=\frac{q_{n_{0}}}{2}$.

    Moving onto the intervals in case (ii), the remaining interval types to handle are the $q_{n_{0}+1,i-1}$-type and the $(q_{n_{0}+1,i-1}-q_{n_{0}})=q_{n_{0}+1,i-2}$-type intervals. We note that these have length at least $\|q_{n_{0}+1,i-1}\alpha\|=(a_{n_{0}+1}-(i-1))\|q_{n_{0}}\alpha\|+\|q_{n_{0}+1}\alpha\|$. By Theorem \ref{Thm:vanRavenstein}, the points $r\alpha$ with $0\leq r <q_{n_{0}+1,i}+q_{n_{0}}$ partition $S^{1}$ into gaps shorter than $\|q_{n_{0}+1,i-1}\alpha\|$, so Lemma \ref{Lem:BoundsOnRho}$(ii)$ gives $\rho(I)\leq \lceil\frac{1}{2}(q_{n_{0}+1,i}+q_{n_{0}})\rceil\leq q_{n_{0}+1,i}$. On the other hand, Lemma \ref{Lem:BoundsOnRho}$(i)$ gives $\rho(I)>N\geq \frac{1}{2}q_{n_{0}+1,i-1}$. We note that $q_{n_{0}+1,i}=q_{n_{0}+1,i-1}+q_{n_{0}}\leq 2q_{n_{0}+1,i-1}$, so with $T=\frac{q_{n_{0}+1,i-1}}{2}$, we have that \eqref{Eqn:T<rho(I)<4T} holds.
   
\end{proof}

\section{Single-level covers and the proof of Theorem \ref{Thm:main}}\label{Sec:UpperIntDimensions}
In this section, we cover $\OmegaAD$ using the complementary intervals of $\J_{N}$ at a single level $N$, and use these covers to prove Theorem \ref{Thm:main}. Subsection \ref{SubSec:LevelStruc} describes the structure of $\J_{N}$ at the levels that matter, Subsection \ref{SubSec:SingleLevel} determines the cheapest such cover at each level, Subsection \ref{SubSec:UB} optimises over the level to obtain the upper bound, and Subsection \ref{SubSec:LB} gives the lower bound.

\subsection{The level structure}\label{SubSec:LevelStruc}

We begin by justifying why the levels of interest lie strictly between the two-gap times 
\begin{equation*}
    N_{n}=\left\lfloor\frac{q_{n}+q_{n+1}}{2} \right\rfloor
\end{equation*}
defined as in \cite{UrbanskiPawalek}. At these times, by Proposition \ref{2GapCase} and Corollary \ref{Cor:Integer3GapThrm}, the complementary intervals of $\J_{N_{n}}$ are $q_{n+1}+O(1)$ intervals of $q_{n}$-type and $q_{n}+O(1)$ intervals of $q_{n+1}$-type, the error terms accounting for parity.

\begin{rmk}\label{rmk:Nn not optimal}
    By Lemma \ref{Lem:BoundsOnRho}, each $q_{n}$-type interval $I$ of $\J_{N_{n}}$ has $\rho(I)>N_{n}\geq \frac{1}{2}q_{n+1}-1$, and $\rho(I)\leq 2q_{n+1}$, so Corollary \ref{cor:BigAppendixCor} gives $|I|_{\dOm}=\Theta\left(q_{n+1}^{-1/\delta}\right)$. Hence, any cover of $\OmegaAD$ in which each complementary interval of $\J_{N_{n}}$ lies in its own covering set has $\sum_{k}|U_{k}|_{\dOm}^{s}\geq c\,q_{n+1}^{1-s/\delta}$, which does not tend to $0$ for $s<\delta$.
    \end{rmk}

From Remark \ref{rmk:Nn not optimal}, we see that this choice of covering time will not illuminate the structure of the minimal set any better than the box dimension was able to. So, if there are any gains to be made, they must occur when considering the system at some intermediate time, $N_{n}<N<N_{n+1}$. With this in mind, we will define the intermediate 2-gap times as follows. We let
\begin{equation*}
    N_{n,i}:= \left\lfloor\frac{q_{n}+(i+1)q_{n+1}}{2}\right\rfloor \;\; \text{ for }\;\; 0\leq i\leq a_{n+2},
\end{equation*}
where we note that $N_{n,0}=N_{n}$, and $N_{n,a_{n+2}}=N_{n+1}$. We will consider times $N$ with $N_{n,i}\leq N <N_{n,i+1}$. Corollary \ref{Cor:Integer3GapThrm} tells us that for $N_{n,i}\leq N <N_{n,i+1}$, the complementary intervals of $\J_{N}$ consist of three types: the $q_{n+1}$-type, $q_{n+2,i}$-type, and the $q_{n+2,i+1}$-type intervals. In the notation from Section \ref{Sec:Approx}, for $i\leq a_{n+1}-2$ the range $N_{n,i}\leq N <N_{n,i+1}$ corresponds, up to parity, to $q_{n_{0}+1,i'-1}<2N+1\leq q_{n_{0}+1,i'}$ with $n_{0}=n+1$ and $i'=i+2$, so the results of that section apply with this substitution. As we run through from $N_{n,i}$ to $N_{n,i+1}$, we split the $q_{n+2,i}$-type intervals, at each step, into a $q_{n+1}$-type interval, and a $q_{n+2,i+1}$-type interval. We let $k=k(N)=2(N-N_{n,i})$ denote the number of $q_{n+2,i}$-type intervals that have been split at time $N$. Then, we can list the number of intervals of each type at time $N$ as follows:
\begin{itemize}
   \item $\#q_{n+2,i}$-type intervals = $q_{n+1}-k+ O(1)$ 
    \item $\#q_{n+2,i+1}$-type intervals = $k+O(1)$
    \item $\#q_{n+1}$-type intervals = $q_{n}+iq_{n+1}+k+O(1)$,
\end{itemize}
where the error terms account for a parity correction.
We note that when $N=N_{n,i}$, we have $k=0$, and when $N=N_{n,i+1}$, we have $k = q_{n+1}$.

We let $\Delta_{n+2,j}$ and $\Gamma_{n+1}$ denote the maximal $\dOm$-lengths of a $q_{n+2,j}$-type interval and a $q_{n+1}$-type interval respectively. By Proposition \ref{Prop:BigAppendixProp}, there is a constant $C \geq 1$, depending only on $\delta$, such that for $1 \leq j \leq a_{n+2}/2$ every $q_{n+2,j}$-type interval has $\dOm$-length at least $C^{-1}\Delta_{n+2,j}$, every $q_{n+1}$-type interval has $\dOm$-length at least $C^{-1}\Gamma_{n+1}$, and
\begin{equation}\label{Eqn:3.1}
  \begin{gathered}
    C^{-1} j^{1-1/\delta}q_{n+1}^{-1/\delta} \leq \Delta_{n+2,j} \leq C j^{1-1/\delta}q_{n+1}^{-1/\delta}, \\
    C^{-1} q_{n+2}^{-1/\delta} \leq \Gamma_{n+1} \leq C q_{n+2}^{-1/\delta}.
  \end{gathered}
\end{equation}

At level $N_{n,i}$, the covering sum will have two competing components whose dependence on $i$ will work against one another: on the one hand, for every $i$, there will always be $q_{n+1}$ $q_{n+2,i}$-type intervals, but their size will be decreasing in $i$ like $i^{1-1/\delta}$; on the other hand, the $q_{n+1}$-type intervals will maintain the same size (roughly $q_{n+2}^{-1/\delta}$), but their number will grow like $iq_{n+1}$. 

The $\theta$-intermediate dimension will limit us to covering $\OmegaAD$ with sets of length in the range $[\Delta^{1/\theta}, \Delta]$ for some $\Delta\in(0,1)$ at any given level. We will use the following lemma to understand the distribution of the intervals making up $\J_{N}$ to help us decide how best to cover them.

\begin{lem}\label{Lem:ClusteringLemma}
        Let $n\geq 1$ be large, let $1\leq i \leq a_{n+2}-1$, and let $N_{n,i}\leq N<N_{n,i+1}$. Then the complementary intervals of $\J_{N}$ admit the following description.
        \begin{enumerate}[(i)]
            \item Each of the $q_{n+1}$ long intervals has a run of consecutive $q_{n+1}$-type intervals on each side of it, with any two consecutive members of a run, and the long interval and each member adjacent to it, being separated by a single Denjoy interval. The two runs lying between consecutive long intervals are separated from one another by a single Denjoy interval $J_{m}$ with $-\frac{q_{n+1}}{2}<m\leq\frac{q_{n+1}}{2}$. Each long interval is of $q_{n+2,i}$-type or of $q_{n+2,i+1}$-type. We call a run, together with the Denjoy intervals strictly between its members, a cluster.
            
            \item Each run has at least $\lfloor\frac{i}{2}\rfloor$ and at most $\lceil\frac{i}{2}\rceil+1$ members.

            \item The Denjoy interval separating the $j^{th}$ and the $(j+1)^{th}$ member of a run (counting towards the long interval), or separating the $j^{th}$ member from the long interval in the run if the run has $j$ members, is $J_{m}$ with $|m|=\Theta\left( jq_{n+1} \right)$, and so has length $\Theta\left(\left(jq_{n+1} \right)^{-1/\delta} \right)$.

           \item Hence, restricting to $2\leq i\leq \frac{a_{n+2}}{2}$, each cluster has $\dOm$-diameter $O\left(q_{n+1}^{-1/\delta}\right)$.
        \end{enumerate}
\end{lem}
\begin{proof}
    We will work in $S^1$ as partitioned by $\{\Ra^{k}(0)\}_{k=-N}^{N}$ and transfer to $\OmegaAD$ by the semiconjugacy, under which each gap becomes a complementary interval of $\J_{N}$ of type determined by the gap length, and the point $m\alpha$ falling between two gaps is blown up to a Denjoy interval $J_{m}$.

    We fix the direction of travel around $S^1$ in which $(m+q_{n+!})\alpha$ lies at distance $\|q_{n+1}\alpha\|$ ahead of $m\alpha$ for every $m$. Then points$\{R^{k}_{\alpha}(0)\}_{k=-N}^{N}$ are the image of $\{\Ra^{k}(0)\}^{2N}_{k=0}$ under the rotation $\Ra^{-N}$, and $q_{n+2,i+1}\leq 2N+1<q_{n+2,i+2}$. So, by Theorems \ref{SteinhausConj} and \ref{Thm:vanRavenstein} applied to $\{\Ra^{k}(0)\}_{k=0}^{2N}$ and read in this direction, for $-N\leq m \leq N$, the next point of the orbit after $m\alpha$ is
    \begin{equation}\label{Cases:ClusterGaps}
        \begin{cases}
            (m+q_{n+1})\alpha & m\leq N-q_{n+1}, \\
            (m-q_{n+2,i+1})\alpha & m>N-q_{n+1} \text{ and }m\geq q_{n+2,i+1}-N, \\ 
            (m-q_{n+2,i})\alpha & N-q_{n+1}<m<q_{n+2,i+1}-N,
        \end{cases}
    \end{equation}
    at distance $\|q_{n+1}\alpha\|$, $\|q_{n+2,i+1}\alpha\|$ and $\|q_{n+2,i}\alpha\|$ respectively. In particular, using the notation from Theorem \ref{SteinhausConj} with $2N+1$ points, Theorem \ref{Thm:vanRavenstein} gives us that $u_{2}=q_{n+1}$ and $u_{2N+1}=q_{n+2,i+1}$ when $2N+1>q_{n+2,i+1}$; and when $2N+1 = q_{n+2,i+1}$, we have $u_{2}=q_{n+1}$ and $u_{2N+1}=q_{n+2,i}$, in which case the second line of \eqref{Cases:ClusterGaps} does not occur since the condition $m\geq q_{n+2,i+1}-N$ would read $m\geq (2N+1)-N=N+1$, and no such index in the orbit satisfies this since $m\in[-N,N]$. The direction fixed above is that of $Suc$ if $u_{2}=q_{n+1}$ and that of $Pre$ otherwise, and \eqref{Cases:ClusterGaps} follow from replacing the index $k$ of a point of $\{\Ra^{k}(0)\}_{k=0}^{2N}$ by $m=k-N$. While the sizes of the gaps are preserved by the rotation $\Ra^{-N}$, the indices of their endpoints, and hence the lengths of the Denjoy intervals between them, are not, which is why we must work with $\{\Ra^{k}(0)\}_{k=-N}^{N}$ directly.

    By \eqref{Cases:ClusterGaps}, travelling in this direction from any point of the orbit, the index increases by $q_{n+1}$ across each gap of length $\|q_{n+1}\alpha\|$ until it exceeds $N-q_{n+1}$. The next gap is then long, and we claim that it leads to an index $m'$ with $m'-q_{n+1}<-N$, so that $m'$ cannot be reached from another point of the orbit across a gap of length $\|q_{n+1}\alpha\|$. There are two cases, according to which one of the last two lines of \eqref{Cases:ClusterGaps} applies. If $m'=m-q_{n+2,i+1}$, then $m'-q_{n+1}=m-q_{n+2,i+2}$, which is less than $-N$ because $m\leq N$ and $2N<q_{n+2,i+2}$. If $m'=m-q_{n+2,i}$, then the second line did not apply, so $m<q_{n+2,i+1}-N$, and hence $m'-q_{n+1}=m-q_{n+2,i+1}<-N$.
    
    So, after each long gap, the indices increase from some $m'<-N+q_{n+1}$ in steps of $q_{n+1}$ until they exceed $N-q_{n+1}$, at which point the next long gap occurs. Between two consecutive long gaps, the indices are therefore all of the integers in $[-N,N]$ which differ from $m'$ by a multiple of $q_{n+1}$, in increasing order, and consecutive points are separated by gaps of length $\|q_{n+1}\alpha\|$. Since $N\geq q_{n+1}$, exactly one of these indices, $m_{0}$ say, satisfies $-\frac{q_{n+1}}{2}<m_{0}\leq\frac{q_{n+1}}{2}$. The gaps before $m_{0}\alpha$ form the run following one long gap, and the gaps after it form the run preceding the next. This proves $(i)$; see Figure \ref{fig:cluster}.
 
    For $(ii)$ and $(iii)$, consider first the run before the long interval. It comes from the points $m_{0},m_{0}+q_{n+1},\dots,m_{0}+kq_{n+1}$, where $N-q_{n+1}<m_{0}+kq_{n+1}\leq N$, so it has $k$ members, and the Denjoy interval following its $j^{th}$ member is $J_{m}$ with $m=m_{0}+jq_{n+1}$. Since $|m_{0}|\leq\frac{q_{n+1}}{2}$, we have that
    \begin{equation*}
        \left(j-\tfrac{1}{2}\right)q_{n+1}\leq |m|\leq \left(j+\tfrac{1}{2}\right)q_{n+1},
    \end{equation*}
    so $|m|=\Theta(jq_{n+1})$ and $J_{m}$ has length $c_{\delta}(1+|m|)^{-1/\delta}=\Theta\left(\left(jq_{n+1}\right)^{-1/\delta}\right)$. Moreover, the bounds on $m_{0}$ and $m_{0}+kq_{n+1}$ give
    \begin{equation*}
        \frac{N}{q_{n+1}}-\frac{3}{2}<k<\frac{N}{q_{n+1}}+\frac{1}{2}.
    \end{equation*}
    Since $N_{n,i}\leq N<N_{n,i+1}$ and $1\leq q_{n}<q_{n+1}$, we have that $\frac{i+1}{2}q_{n+1}\leq N<\frac{i+3}{2}q_{n+1}$. Hence $\frac{i}{2}-1<k<\frac{i}{2}+2$, giving us that, $\left\lfloor\frac{i}{2}\right\rfloor\leq k\leq\left\lceil\frac{i}{2}\right\rceil+1$. The run after the long interval comes from the points $m_{0}-kq_{n+1},\ldots,m_{0}-q_{n+1},m_{0}$ of the next chain, where $-N\leq m_{0}-kq_{n+1}<-N+q_{n+1}$, and the same argument applies with $m=m_{0}-jq_{n+1}$. This proves $(ii)$ and $(iii)$.

    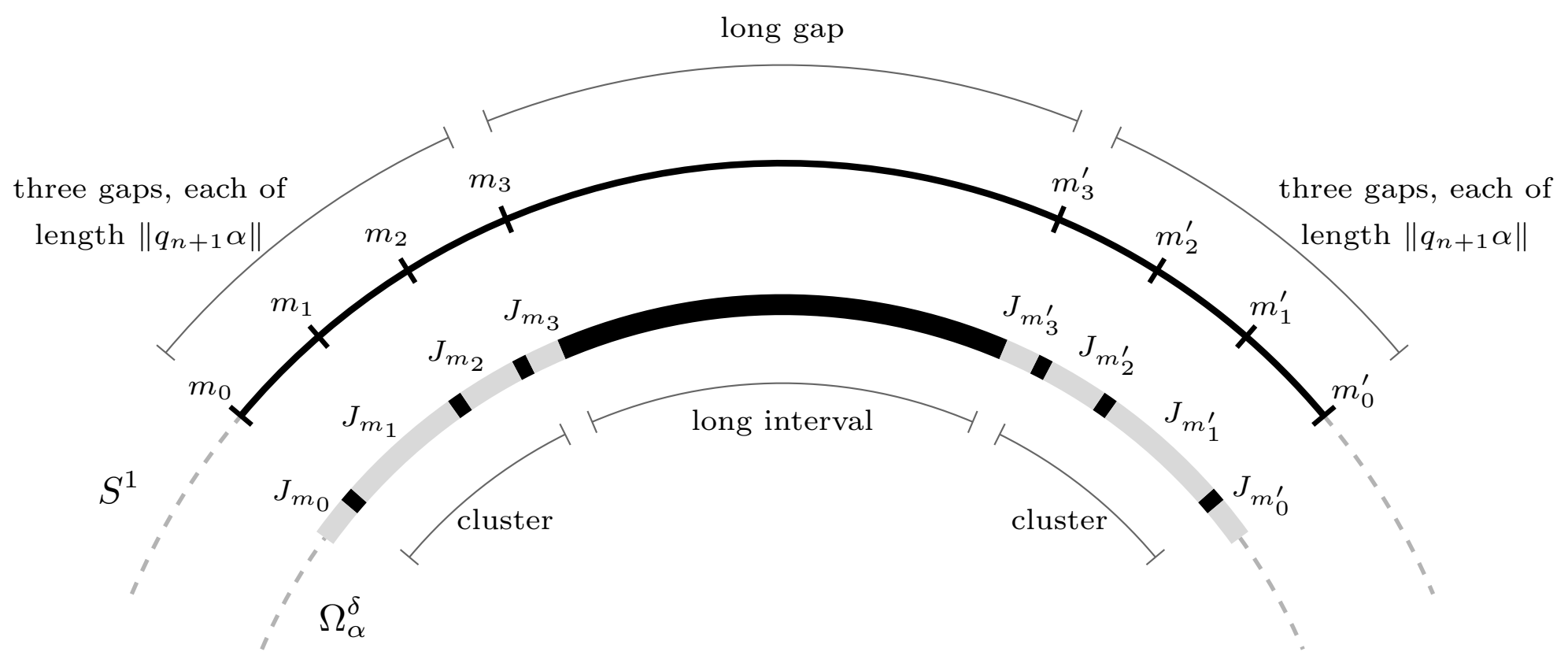
\begin{figure}[ht]
\centering
\begin{tikzpicture}[scale=0.66]

\def\Rone{9}    
\def\Rtwo{7.2}  

\draw[black!30,dashed,line width=1pt] (157:\Rone) arc (157:140:\Rone);
\draw[black!30,dashed,line width=1pt] (40:\Rone)  arc (40:23:\Rone);
\draw[black!30,dashed,line width=1pt] (157:\Rtwo) arc (157:144:\Rtwo);
\draw[black!30,dashed,line width=1pt] (36:\Rtwo)  arc (36:23:\Rtwo);

\draw[line width=1.6pt] (140:\Rone) arc (140:40:\Rone);

\foreach \a in {140,131,122,113,67,58,49,40}
  \draw[line width=1.2pt] (\a:\Rone-0.18) -- (\a:\Rone+0.18);

\node at (140:\Rone+0.5) {\scriptsize $m_{0}$};
\node at (131:\Rone+0.5) {\scriptsize $m_{1}$};
\node at (122:\Rone+0.5) {\scriptsize $m_{2}$};
\node at (113:\Rone+0.5) {\scriptsize $m_{3}$};
\node at (67 :\Rone+0.5) {\scriptsize $m_{3}'$};
\node at (58 :\Rone+0.5) {\scriptsize $m_{2}'$};
\node at (49 :\Rone+0.5) {\scriptsize $m_{1}'$};
\node at (40 :\Rone+0.5) {\scriptsize $m_{0}'$};

\foreach \a/\b in {140/114.5, 111.5/68.5, 65.5/40}{
  \draw[black!60] (\a:\Rone+1.25) arc (\a:\b:\Rone+1.25);
  \draw[black!60] (\a:\Rone+1.1) -- (\a:\Rone+1.4);
  \draw[black!60] (\b:\Rone+1.1) -- (\b:\Rone+1.4);
}

\node[align=center] at (134:\Rone+2.6)
  {\scriptsize three gaps, each of\\[-1pt] \scriptsize length $\|q_{n+1}\alpha\|$};
\node at (90:\Rone+1.7) {\scriptsize long gap};
\node[align=center] at (46:\Rone+2.6)
  {\scriptsize three gaps, each of\\[-1pt] \scriptsize length $\|q_{n+1}\alpha\|$};

\node at (150:\Rone+0.75) {\small $S^{1}$};

\foreach \a/\b in {140/138.5, 125.5/124, 118/116.5, 113/67, 63.5/62, 56/54.5, 41.5/40}
  \draw[line width=5pt] (\a:\Rtwo) arc (\a:\b:\Rtwo);

\foreach \a/\b in {138.5/125.5, 124/118, 116.5/113, 67/63.5, 62/56, 54.5/41.5}
  \draw[line width=5pt,black!15] (\a:\Rtwo) arc (\a:\b:\Rtwo);

\draw[line width=5pt,black!15] (144:\Rtwo) arc (144:140:\Rtwo);
\draw[line width=5pt,black!15] (40:\Rtwo)  arc (40:36:\Rtwo);
\node at (142:\Rtwo+0.55) {\scriptsize $J_{m_{0}}$};
\node at (38 :\Rtwo+0.55) {\scriptsize $J_{m_{0}'}$};

\node at (132.5:\Rtwo+0.55) {\scriptsize $J_{m_{1}}$};
\node at (122.3:\Rtwo+0.55) {\scriptsize $J_{m_{2}}$};
\node at (114.2:\Rtwo+0.55) {\scriptsize $J_{m_{3}}$};
\node at (65.8:\Rtwo+0.55) {\scriptsize $J_{m_{3}'}$};
\node at (57.7:\Rtwo+0.55) {\scriptsize $J_{m_{2}'}$};
\node at (47.5:\Rtwo+0.55) {\scriptsize $J_{m_{1}'}$};

\draw[black!60] (140:\Rtwo-1.0) arc (140:116.5:\Rtwo-1.0);
\draw[black!60] (140  :\Rtwo-0.85) -- (140  :\Rtwo-1.15);
\draw[black!60] (116.5:\Rtwo-0.85) -- (116.5:\Rtwo-1.15);
\node at (128.25:\Rtwo-1.5) {\scriptsize cluster};

\draw[black!60] (113:\Rtwo-1.0) arc (113:67:\Rtwo-1.0);
\draw[black!60] (113:\Rtwo-0.85) -- (113:\Rtwo-1.15);
\draw[black!60] (67 :\Rtwo-0.85) -- (67 :\Rtwo-1.15);
\node at (90:\Rtwo-1.5) {\scriptsize long interval};

\draw[black!60] (63.5:\Rtwo-1.0) arc (63.5:40:\Rtwo-1.0);
\draw[black!60] (63.5:\Rtwo-0.85) -- (63.5:\Rtwo-1.15);
\draw[black!60] (40  :\Rtwo-0.85) -- (40  :\Rtwo-1.15);
\node at (51.75:\Rtwo-1.5) {\scriptsize cluster};

\node at (150:\Rtwo-0.75) {\small $\Omega^{\delta}_{\alpha}$};

\end{tikzpicture}

\caption{A long interval and the two runs adjacent to it, drawn with three members in each run. Outer arc: the points of the orbit in $S^{1}$, with indices $m_{j}=m_{0}+jq_{n+1}$ and $m_{j}'=m_{0}'-jq_{n+1}$. Inner arc: the corresponding intervals of $\Omega^{\delta}_{\alpha}$, with the Denjoy intervals shaded. The two arcs are drawn to different scales.}
\label{fig:cluster}
\end{figure}

    Finally, for $(iv)$, we first note that the $\dOm$-length of a cluster of $q_{n+1}$ type intervals will be at most
    \begin{equation*}
        i\Theta\left(q_{n+2}^{-1/\delta} \right) 
        + \Theta\left(\sum_{j=1}^{i-1} (jq_{n+1})^{-1/\delta}\right).
    \end{equation*}
    Here, the first term comes from the lengths of the at most $i$ intervals of $q_{n+1}$-type, and the second comes from the lengths of the Denjoy intervals separating them. Analysing the second term, we note that
    \begin{align*}
        \sum_{j=1}^{i-1} (jq_{n+1})^{-1/\delta} &=q_{n+1}^{-1/\delta}\sum_{j=1}^{i-1}j^{-1/\delta} \\ &\leq q_{n+1}^{-1/\delta}\sum_{j=1}^{\infty}j^{-1/\delta} \\ &= \zeta(1/\delta)q_{n+1}^{-1/\delta}.
    \end{align*}
    Next, we note that since $q_{n+2}=\Theta(a_{n+2}q_{n+1})$, the $q_{n+1}$-type intervals will contribute
    \begin{align*}
        i\Theta\left(q_{n+2}^{-1/\delta} \right) &= \Theta\left(iq_{n+2}^{-1/\delta} \right) \\ &= \Theta\left(\frac{i}{(a_{n+2}q_{n+1})^{1/\delta}} \right) \\ &=\Theta\left(\frac{i}{a_{n+2}^{1/\delta}}q_{n+1}^{-1/\delta} \right),
    \end{align*}
    so since $i\leq a_{n+2}$, we have that $\frac{i}{a_{n+2}^{1/\delta}}\leq a_{n+2}^{1-1/\delta}$. Then, since $1/\delta >1$, $1-1/\delta<0$, and so the contribution of the $q_{n+1}$-type intervals to the overall length of the cluster is $O\left(q_{n+1}^{-1/\delta} \right)$. Together with the bound on the Denjoy intervals above, this proves $(iv)$.
    \end{proof}

    \subsection{Single-level covers}\label{SubSec:SingleLevel}
    
    Before moving on to considering the lower $\theta$-intermediate dimension, we recall the following notation defined earlier:
    \begin{equation*}
    \Delta_{n+2,i}:=\max\left\{|I|_{\dOm}:I \text{ is of } q_{n+2,i}\text{-type} \right\}    \text{ and } \Gamma_{n+1}:=\max\left\{|I|_{\dOm}:I \text{ is of } q_{n+1}\text{-type} \right\}.
    \end{equation*}
    Defining these lengths allows us to ensure any covering set for a $q_{n+2,i}$ or $q_{n+1}$-type interval will be covered by a set of size $\Delta_{n+2,i}$ or $\Gamma_{n+1}$, respectively, and since all such sets are of comparable size, using such covering sets will not inflate the dimension estimates.
    \begin{definition}\label{Def:LevelNCover}
        Let $\theta\in (0,1]$, let $n\geq 2$ and $i\geq 2$ with $i+1\leq \frac{a_{n+2}}{2}$, and let $N\in\mathbb{N}$ satisfy $N_{n,i}\leq N<N_{n,i+1}$. Put $\Delta=\Delta_{n+2,i}$, and $\beta:=\max\left\{\Gamma_{n+1},\Delta^{1/\theta} \right\}$. A cover $\{U_{k}\}$ of $\OmegaAD$ is a \emph{level-}$N$ \emph{cover for} $\theta$ if 
        \begin{enumerate}
            \item $\beta\leq |U_{k}|_{\dOm}\leq \Delta$ for every $k$, and
            \item every connected component of $\J_{N}$ is contained in some $U_{k}$.
        \end{enumerate}
    \end{definition}

    \begin{rmk}\label{Rmk:BetaDelta}
        Condition $2$ above says no complementary interval is split between two covering sets; it permits a single $U_{k}$ to contain several complementary intervals, which is the strategy tested below.

        We also note that level-$N$ covers for $\theta$ only exist if $\beta\leq \Delta$. Since $\Delta<1$, we always have that $\Delta^{1/\theta}\leq\Delta$, and $\Gamma_{n+1}\leq \Delta$ holds once $a_{n+2}$ is large enough in terms of $\delta$. Indeed, by Proposition \ref{Prop:BigAppendixProp} and the fact that $a_{n+2}q_{n+1}<q_{n+2}$, 
        \begin{equation*}
            \Gamma_{n+1}\leq C(a_{n+2}q_{n+1})^{-1/\delta} \;\; \text{ and } \;\; \Delta\geq ci^{1-1/\delta}q_{n+1}^{-1/\delta}  \geq c\left(\frac{a_{n+2}}{2}\right)^{1-1/\delta}q_{n+1}^{-1/\delta},
        \end{equation*}
        giving us that $\frac{\Delta}{\Gamma_{n+1}}\geq ca_{n+2}$. Here, $C,c>0$ depend only on $\delta$.
    \end{rmk}
   Before proving the main theorem of this section, we pause to determine the most efficient way to cover the $q_{n+1}$-type intervals of $\J_{N}$.

\begin{lem}\label{Lem:SingleLevelCost}
Let $\alpha \in (0,1)$ be irrational, let $\delta, s \in (0,1)$, and let $\theta$, $n$, $i$, $N$, $\Delta$ and $\beta$ be as in Definition~\ref{Def:LevelNCover}. Suppose that $a_{n+2}$ is sufficiently large in terms of $\delta$. Then
\begin{equation*}
  \inf \sum_{k} |U_{k}|_{\dOm}^{s} = \Theta\left( q_{n+1}\Delta^{s} +  \min\left\{iq_{n+1}\beta^{s}, \beta^{s(1-\delta)}\Delta^{-\delta(1-s)}  \right\} \right),
\end{equation*}
where the infimum is taken over all level-$N$ covers $\{U_{k}\}$ for $\theta$, and the implied constants depend only on $\delta$.
\end{lem}
 
\begin{rmk}
The three terms correspond to the three parts of an optimal cover. The term $q_{n+1}\Delta^{s}$ is the cost of the long intervals. The term $iq_{n+1}\beta^{s}$ is the cost of covering every $q_{n+1}$-type interval individually by its own set of diameter $\beta$, and the term $\beta^{s(1-\delta)}\Delta^{-\delta(1-s)}$ is the cost of sets of diameter $\Delta$ when shared between $q_{n+1}$-type intervals in the case where this is cheaper.
\end{rmk}
 
\begin{proof}
Throughout, $c$ and $C$ denote positive constants depending only on $\delta$, whose values may change from line to line. By Remark \ref{Rmk:BetaDelta} we have $\beta \leq \Delta$. Replacing a covering set by the smallest arc containing it does not increase its $\dOm$-diameter, so we may take every covering set to be an arc; for the same reason, an arc of diameter at most $\Delta$ may be enlarged to one of diameter in $[\beta,\Delta]$.
 
To begin with, we establish the lengths of intervals we need to cover. By Lemma \ref{Lem:ClusteringLemma}, the $q_{n+1}$-type intervals of $\J_{N}$ fall into $2q_{n+1}$ clusters, one on each side of each long interval, which is of $q_{n+2,i}$-type or $q_{n+2,i+1}$-type. Since $i\geq 2$, each cluster has at least $\lfloor\frac{i}{2}\rfloor\geq 1$ and at most $\lceil\frac{i}{2}\rceil +1\leq i$ member, which we number towards the long interval. We let $g_{j}$ denote the $\dOm$-length of the Denjoy interval following the $j^{th}$ member of a cluster, and have that
\begin{equation*}
  g_{j} = \Theta\big((jq_{n+1})^{-1/\delta}\big), \;\;\;\;\; 1 \leq j \leq i.
\end{equation*}
By Proposition \ref{Prop:BigAppendixProp}, since $i+1 \leq a_{n+2}/2$, the $q_{n+2,i}$-type and $q_{n+2,i+1}$-type intervals have $\dOm$-lengths $\Theta\big(i^{1-1/\delta}q_{n+1}^{-1/\delta}\big)$ and $\Theta\big((i+1)^{1-1/\delta}q_{n+1}^{-1/\delta}\big)$ respectively. As $(i+1)^{1-1/\delta} \geq 2^{1-1/\delta}i^{1-1/\delta}$, every long interval has $\dOm$-length at least $c\Delta$, and
\begin{equation}\label{eq:Delta}
  \Delta = \Theta\big(i^{1-1/\delta}q_{n+1}^{-1/\delta}\big) = \Theta(i\,g_{i}).
\end{equation}
 
There are $q_{n+1}$ long intervals. An arc of diameter at most $\Delta$ containing $k$ of them has diameter at least $kc\Delta$, so it contains at most $1/c$ of them; and since consecutive clusters are separated alternately by a long interval and by a single Denjoy interval, an arc of diameter at most $\Delta$ meeting $2k+1$ clusters contains at least $k$ long intervals entirely, so it meets at most $2/c + 2$ clusters. In a level-$N$ cover for $\theta$, each long interval lies in a set of diameter at least $c\Delta$, and each set contains at most $1/c$ of them, so
\begin{equation}\label{eq:long-lower}
  \sum_{k} |U_{k}|_{\dOm}^{s} \;\geq\; c\,q_{n+1}\Delta^{s}.
\end{equation}
Conversely, giving each long interval its own set, of diameter at most $\Delta$, costs at most $q_{n+1}\Delta^{s}$.
 
We now make two observations. First, the Denjoy intervals inside a single covering set have comparable lengths. Let $U$ be an arc of diameter at most $\Delta$ meeting a single cluster, and let $g_{j},\dots,g_{m}$ be the lengths of the Denjoy intervals of that cluster contained in $U$. Since $g_{l} \geq c\,(lq_{n+1})^{-1/\delta}$, comparing with an integral gives
\begin{align*}
  \Delta &\geq \sum_{l=j}^{m} g_{l} \\&\geq \int_{j}^{m+1} c(xq_{n+1})^{-1/\delta}\,dx \\&=\; c'q_{n+1}^{-1/\delta}\left(j^{1-1/\delta}-(m+1)^{1-1/\delta}\right),
\end{align*}
where $c' = c/(1/\delta - 1)$. By \eqref{eq:Delta}, $\Delta \leq C i^{1-1/\delta}q_{n+1}^{-1/\delta}$. Put $\lambda = \big(c'/(2C)\big)^{\delta/(1-\delta)}$.
\begin{itemize}
  \item If $j \leq \lambda i$, then $C i^{1-1/\delta} \leq \tfrac{c'}{2} j^{1-1/\delta}$, so $(m+1)^{1-1/\delta} \geq \tfrac{1}{2}j^{1-1/\delta}$, that is, $m+1 \leq 2^{\delta/(1-\delta)}j$.
  \item If $j > \lambda i$, then, since $U$ meets a single cluster, $m \leq i< 2j/\lambda$.
\end{itemize}
So $m \leq Kj$ with $K = \max\{2^{\delta/(1-\delta)}, 2/\lambda\}$, which depends only on $\delta$. Since we have that $g_{l} = \Theta\big((lq_{n+1})^{-1/\delta}\big)$, it follows that
\begin{equation*}
  c\,g_{j} \leq g_{l} \leq C g_{j} \;\; \text{ whenever } j \leq l \leq Kj ,
\end{equation*}
and in particular for every Denjoy interval in $U$. We also note that $g_{l} \leq C g_{j}$ for every $l \geq j$.
 
Second, the members of a cluster are no longer than the Denjoy intervals between them: for every $l \leq a_{n+2}$, by Proposition \ref{Prop:BigAppendixProp} and $a_{n+2}q_{n+1} \leq q_{n+2}$,
\begin{equation*}
  \Gamma_{n+1} \;\leq\; C(a_{n+2}q_{n+1})^{-1/\delta} \;\leq\; C(lq_{n+1})^{-1/\delta} \;\leq\; C g_{l}.
\end{equation*}
 
An arc meets the members of a cluster in a run of consecutive members, and we say it \emph{starts at} the $j^{th}$ member if that is the first member it contains. The two observations above tell us that within a covering set, every member and every Denjoy interval has length at most $Cg_{j}$, and every Denjoy interval has length at least $cg_{j}$. So, if $L$ is the $\dOm$-length of the smallest arc containing the $j^{th}$ and the $(j+m-1)^{th}$ members, which consist of these $m$ members and the $m-1$ Denjoy intervals between them, then 
\begin{equation*}
    (m-1)cg_{j}\leq L\leq mCg_{j},
\end{equation*}
 where the lower bound requires $L\leq \Delta$. This gives the number of members that a set of diameter $D\leq \Delta$ can hold. On the one hand, an arc of diameter $D$ starting at the $j^{th}$ member containing $m$ members must have that $D\geq L$, so by the equation above,
 \begin{equation*}
     m\leq 1+\frac{D}{cg_{j}}.
 \end{equation*}
On the other hand, if $D\geq \Gamma_{n+1}$ and the cluster has at least
\begin{equation}\label{Eqn:MembersInDiamDCoverSet}
    m=\max\left\{ 1,\left\lfloor \frac{D}{Cg_{j}}\right\rfloor\right\} 
\end{equation}
 members starting from the $j^{th}$, then these $m$ members lie in an arc of diameter at most $D$. To see this, we note that if $m=1$, then a single member has length at most $\Gamma_{n+1}\leq D$. When $m\geq 2$, we have that $m=\left\lfloor\frac{D}{Cg_{j}} \right\rfloor$, meaning that $m\leq \frac{D}{Cg_{j}}$ and hence $D\geq mCg_{j}\geq L$, and thus must contain the $j^{th}$ to $(j+m-1)^{th}$ members. In other words, a set of diameter $D\in[\Gamma_{n+1},\Delta]$ starting from the $j^{th}$ member holds $\Theta(\max\{1,D/g_{j}\})$ members, with the upper bound always holding, and the lower bound depending on the cluster not running out of members first.

 With this in hand, we now calculate the cost of each $q_{n+1}$-type interval in a cluster to our covering sum. The number of members contained in a set of diameter $D\in[\beta,\Delta]$ starting at the $j^{th}$ member is $\Theta(\max\{1,D/g_{j}\})$, so each such member will take up $\Theta\left( \frac{1}{\max\{1,D/g_{j}\}}\right)$ covering sets, each of which costs $D^{s}$ to our covering sum. So, the cost function of the $j^{th}$ member of our covering sum will be 
 \begin{equation*}
     C_{j}(D):=\frac{D^{s}}{\max\{1,D/g_{j}\}} =\min\{D^{s},D^{s-1}g_{j}\}.
 \end{equation*}
We want to minimise this function, i.e. find $\min_{\beta\leq D\leq \Delta}C_{j}(D)$. To do this, we note that the function $D^{s}$ is increasing in $D$, while $D^{s-1}g_{j}$ is decreasing in $D$ since $s<1$. Hence, we can write
\begin{align*}
    \min_{\beta\leq D\leq \Delta}C_{j}(D)&= \min_{\beta\leq D\leq \Delta}\min\{D^{s},D^{s-1}g_{j}\} \\ &= \min \left\{\min_{\beta\leq D\leq \Delta} D^{s},\min_{\beta\leq D\leq \Delta} D^{s-1}g_{j} \right\} \\ &= \min \left\{\beta^{s},\Delta^{s-1}g_{j} \right\}  
\end{align*}

We now want to come up with a lower bound for the covering sum. To begin, let $\{U_{k}\}$ be a level-$N$ cover for $\theta$. By condition $(2)$ from the definition, each member of the cluster must be contained in some covering set. Assign each member to a set accordingly. For each set $U_{k}$ in our cover, we share $\left|U_{k}\right|^{s}_{\dOm}$ among the members of each cluster assigned to $U_{k}$, first equally between the at most $C$ clusters that $U_{k}$ meets, and then equally among the members of each. If $U_{k}$ has diameter $D$ and starts at the $j^{th}$ member of a cluster, then our earlier work tells us that it contains at most $C\max\{1,D/g_{j}\}$ members of the cluster. Each such member, say the $l^{th}$, has $l\geq j$ and so $g_{l}\leq Cg_{j}$. Each member therefore costs at least
\begin{align*}
    c\frac{D^{s}}{\max\{1,D/g_{j}\}} &\geq c\frac{D^{s}}{\max\{1,D/g_{l}\}} \geq c\min\{\beta^{s},\Delta^{s-1}g_{l}\}
\end{align*}
of the mass of $U_{k}$. The costs of the members assigned to $U_{k}$ sum to at most the contribution of $|U_{k}|^{s}_{\dOm}$, and each member belongs to exactly one set in the cover, so the covering sum is bounded below by the total cost of all of the members of each cluster. Since there are $2q_{n+1}$ clusters, each with at least $\lfloor\frac{i}{2}\rfloor$ members, this gives us
\begin{equation*}
    \sum_{k}|U_{k}|^{s}_{\dOm} \geq cq_{n+1}\sum_{j=1}^{\lfloor\frac{i}{2}\rfloor}\min\{\beta^{s},\Delta^{s-1}g_{j}\}.
\end{equation*}
Combining this with the covering cost for the larger intervals, and using the fact that $\max\{A,B\}\geq \frac{1}{2}(A+B)$, we get that
\begin{equation}\label{Eqn:CovSumLB}
    \sum_{k}|U_{k}|^{s}_{\dOm} \geq cq_{n+1}\left(\Delta^{s}+\sum_{j=1}^{\lfloor\frac{i}{2}\rfloor}\min\{\beta^{s},\Delta^{s-1}g_{j}\} \right).
\end{equation}
We note that $g_{j}=\Theta\left((jq_{n+1})^{-1/\delta}\right)$ is decreasing in $j$, and so in particular, we want to understand the values of $j$ for which the minimum $\min\{\beta^{s},\Delta^{s-1}g_{j}\}$ takes each value. As such, we want to know when
\begin{equation*}
    \beta^{s}=\Delta^{s-1}g_{j} \iff g_{j}=\beta^{s}\Delta^{1-s}.
\end{equation*}
We write 

\begin{equation*}
    (j_{0}q_{n+1})^{-1/\delta}=\beta^{s}\Delta^{1-s}, \;\; \text{ i.e. }\;\;  j_{0}=q_{n+1}^{-1}\left(\beta^{s}\Delta^{1-s}\right)^{-\delta},
\end{equation*}
a positive real number. Since $g_{j}=\Theta\left((jq_{n+1})^{-1/\delta}\right)$, we have $\Delta^{s-1}g_{j}\geq c\Delta^{s-1}(j_{0}q_{n+1})^{-1/\delta}= c\beta^{s}$ for $j\leq j_{0}$, and that $\Delta^{s-1}g_{j}\leq C\Delta^{s-1}(j_{0}q_{n+1})^{-1/\delta} =C\beta^{s}$ for $j\geq j_{0}$. Thus, the cost of the $j^{th}$ member of a cluster to our covering sum will be $\Theta(\beta^{s})$ for $j\leq j_{0}$, and $\Theta(\Delta^{s-1}g_{j})$ for $j>j_{0}$. We note that the index at which $\Delta^{s-1}g_{j}$ actually passes $\beta^{s}$ may differ from $j_{0}$ by a bounded factor, but we only use the above inequalities in what follows.

We now want to construct a level-$N$ cover for $\theta$ to produce an upper bound on our covering sum. Aside from the $O(q_{n+1})$ larger intervals of diameter at most $\Delta$, the cost of the $j^{th}$ member of each cluster will be at most $C\min\{\beta^{s},\Delta^{s-1}g_{j}\}$. By the discussion introducing $j_{0}$, this means covering the $j^{th}$ member individually when $j\leq j_{0}$, and sharing a set of diameter $\Delta$ among multiple members when $j>j_{0}$. The cover will then consist of the following sets.

\begin{enumerate}[(a)]
    \item \emph{Long intervals}: each long interval will be covered individually by a set of diameter at most $\Delta$, contributing a combined cost of $q_{n+1}\Delta^{s}$.
    \item \emph{Members with $j\leq \min\{i,j_{0}\}$}: each such interval is covered individually by a set of diameter $\beta$, each contributing $\beta^{s}$ to the covering sum. By the argument given to define $j_{0}$, we have that $\beta^{s}\leq C\Delta^{s-1}g_{j}$ for $j\leq j_{0}$, so this cost can be upper bounded by $C\min\{\beta^{s},\Delta^{s-1}g_{j}\}$.
    \item \emph{Members with $j_{0}<j\leq i$}: in each cluster, these are covered in order of increasing index by sets of diameter $\Delta$. The first set starts at the first such member, i.e. the first $j$ for which $j>j_{0}$. If a set starts at the $j^{th}$ member, it contains $m_{j}=\max\{1,\lfloor\Delta/Cg_{j}\rfloor\}$ members from the $j^{th}$ onwards, and the next set starts at the following member. When fewer than $m_{j}$ members remain, the last set instead contains all of them. By our earlier argument, the length of the smallest arc containing $m$ members is bounded above by $\Delta$, and thus so will be the length of the covering set containing all remaining members here. There is at most one such set in each cluster, so these cost at most $2q_{n+1}\Delta^{s}$ in total. We will call a set containing exactly $m_{j}$ members \emph{full} and charge its contribution of $\Delta^{s}$ equally to its $m_{j}$ constituent members. We now claim that $\Delta^{s}\leq 2Cm_{j}\Delta^{s-1}g_{j}$. Indeed, if $\Delta\geq Cg_{j}$, then $m_{j}\geq \frac{\Delta}{2Cg_{j}}$, while if $\Delta<Cg_{j}$, then $m_{j}=1$, and $\Delta^{s}\leq C\Delta^{s-1}g_{j}$. Each member of a full set, say the $l^{th}$, satisfies $j\leq l\leq Kj$, and so $g_{l}\geq cg_{j}$. Each such set therefore costs 
    \begin{equation*}
        \frac{\Delta^{s}}{m_{j}}\leq 2C\Delta^{s-1}g_{j}\leq C\Delta^{s-1}g_{l}\leq C\min\{\beta^{s},\Delta^{s-1}g_{l}\}
    \end{equation*}
    where the last inequality uses the fact that we are in the case where $l\geq j>j_{0}$.
\end{enumerate}
Every set in this proposed cover has diameter in $[\beta,\Delta]$, and every connected component of $\J_{N}$ lies in one of them, so this cover is indeed a level-$N$ cover for $\theta$. The sets in $(a)$, and the sets in $(c)$ which are not full, contribute at most $3q_{n+1}\Delta^{s}$ to the covering sum. By $(b)$ and $(c)$, for $1\leq j\leq i$, the cost of the $j^{th}$ member of each cluster that has one will be at most $C\min\{\beta^{s},\Delta^{s-1}g_{j}\}$. Summing over the $2q_{n+1}$ clusters gives us the following upper bound for the covering sum
\begin{equation*}
    \sum_{k}|U_{k}|^{s}_{\dOm}\leq Cq_{n+1}\left(\Delta^{s}+\sum_{j=1}^{i}
    \min \{\beta^{s},\Delta^{s-1}g_{j}\}\right).
\end{equation*}

We now use the definition of $j_{0}$ to evaluate the summation from the upper and lower bounds. We will show that
\begin{equation}\label{Eqn:EvalSumPreambleLem}
    \Delta^{s}+\sum_{j=1}^{i}\min\{\beta^{s},\Delta^{s-1}g_{j}\} = \Theta(\Delta^{s}+\min\{i,j_{0}\}\beta^{s})
\end{equation}
by proving upper and lower bounds.
By our definition of $j_{0}$, the $j^{th}$ member of a cluster for whom $j\leq\min\{i,j_{0}\}$ will contribute 
\begin{equation*}
    \Theta(\lfloor\min\{i,j_{0}\}\rfloor\beta^{s})
\end{equation*}
to the covering sum. 

We now prove the upper and lower bounds for the remaining members separately. We note that terms with $\min\{i,j_{0}\}< j \leq i$ only exist if $\min\{i,j_{0}\}=j_{0}$, in which case they are terms with $j_{0}< j \leq i$, each at most $\Delta^{s-1}g_{j}$.  If $j_{0}\geq 1$, then 
\begin{align*}
    \sum_{j_{0}\leq j\leq i}\Delta^{s-1}g_{j} &\leq \sum_{j_{0}\leq j\leq i}C\Delta^{s-1}(jq_{n+1})^{-1/\delta}\\ &=C\Delta^{s-1}q_{n+1}^{-1/\delta}\sum_{j_{0}\leq j\leq i}j^{-1/\delta}.
\end{align*}
We note that using an integral comparison, we have that
\begin{align*}
    \sum_{j_{0}\leq j\leq i}j^{-1/\delta} &\leq \sum_{j>j_{0}}j^{-1/\delta}\\ &\leq \int _{\lfloor j_{0}\rfloor}^{\infty}x^{-1/\delta}\;dx \\ &=\frac{\lfloor j_{0}\rfloor^{1-1/\delta}}{1/\delta-1} \leq Cj_{0}^{1-1/\delta}
\end{align*}
where we use that $\lfloor j_{0}\rfloor\geq j_{0}/2$ for $j_{0}\geq 1$ for the last inequality.
So putting this back into our previous string of inequalities, we have that 
\begin{align*}
    \sum_{j_{0}\leq j\leq i}\Delta^{s-1}g_{j}\ &\leq C\Delta^{s-1}q_{n+1}^{-1/\delta}j_{0}^{1-1/\delta} \\ &= C\Delta^{s-1}j_{0}(j_{0}q_{n+1})^{-1/\delta} \\ &= Cj_{0}\beta^{s}
\end{align*}
where the last equality comes from  $(j_{0}q_{n+1})^{-1/\delta}=\beta^{s}\Delta^{1-s}$. We note that, here, $j_{0}=\min\{i,j_{0}\}$, and so these terms contribute at most $C\min\{i,j_{0}\}\beta^{s}$. If instead we have that $j_{0}<1$, then we note that $q_{n+1}^{-1/\delta}<\beta^{s}\Delta^{1-s}$ by the definition of $j_{0}$. So, summing the whole series gives us
\begin{align*}
    \sum_{j=1}^{i}\Delta^{s-1}g_{j}&\leq\Delta^{s-1}\sum_{j=1}^{i}C(jq_{n+1})^{-1/\delta} \\&=C\Delta^{s-1}q_{n+1}^{-1/\delta}\sum_{j=1}^{i}j^{-1/\delta} \\ &\leq C\Delta^{s-1}q_{n+1}^{-1/\delta} \\ &<C\beta^{s} \\&\leq C\Delta^{s}.
\end{align*}
In each case, we have shown that the left-hand side of \eqref{Eqn:EvalSumPreambleLem} is at most $C\left(\Delta^s+\min\{i,j_{0}\}\beta^{s}\right)$.

For the lower bound, we note that by the definition of $j_{0}$, each term in the sum with $j<\min\{\lfloor\frac{i}{2}\rfloor,j_{0}\}$ will be at least $c\beta^{s}$, and there will be $\lfloor\min\{\lfloor\frac{i}{2}\rfloor,j_{0}\}\rfloor$ of them. If we consider the $j=1$ term, since $i\geq 2$ and $\Delta=\Theta(ig_{i})$, we have that $\Delta\leq Ci^{1-1/\delta}q_{n+1}^{-1/\delta} \leq Cq_{n+1}^{-1/\delta} \leq Cg_{1}$. Hence, $\Delta^{s-1}g_{1}\geq \Delta^{s}\geq \beta^{s}$. So, at least $\max\{1,\lfloor\min\{\lfloor\frac{i}{2}\rfloor,j_{0}\}\rfloor\}$ terms take a value of at least $c\beta^{s}$. Since $\max\{1,\lfloor x\rfloor\}\geq x/2$ for any $x\geq0$, and $\lfloor\frac{i}{2}\rfloor\geq \frac{i}{3}$, this gives us that
\begin{equation*}
    \sum_{j=1}^{\lfloor\frac{i}{2}\rfloor}\min\{\beta^{s},\Delta^{s-1}g_{j}\}\geq c\min\{i,j_{0}\}\beta^{s}.
\end{equation*}
The sum in \eqref{Eqn:EvalSumPreambleLem} is at least as large, so adding the leading $\Delta^s$ term to both sides completes the proof of \eqref{Eqn:EvalSumPreambleLem}.

We have shown that any covering sum derived from a level-$N$ cover for $\theta$ will have a form bounded below by the right-hand side of \eqref{Eqn:CovSumLB}, which by the inequality above is at least $cq_{n+1}\left(\Delta^{s}+\min\{i,j_{0}\}\beta^{s} \right)$, and then found a valid level-$N$ cover for $\theta$ whose covering sum has the same form as an upper bound. This is sufficient to say that this form describes an infimum.

Finally, returning to our covering sum and combining it with \eqref{Eqn:EvalSumPreambleLem} we get that
\begin{equation*}
    \inf \sum_{k}|U_{k}|^{s}_{\dOm}=\Theta\left(q_{n+1}\left(  \Delta^{s}+\min\{i,j_{0}\}\beta^{s} \right)\right) = \Theta\left(q_{n+1}\Delta^{s}+q_{n+1}\min\{i,j_{0}\}\beta^{s} \right).
\end{equation*}
The second term in this expression takes the value $\min\left\{iq_{n+1}\beta^{s},j_{0}q_{n+1}\beta^{s}\right\}$. The latter of these terms can be simplified by noting that $j_{0}q_{n+1}=\left(\beta^{s}\Delta^{1-s} \right)^{-\delta}$, which gives us that
\begin{equation*}
    j_{0}q_{n+1}\beta^{s}=\left(\beta^{s}\Delta^{1-s} \right)^{-\delta}\beta^{s}=\beta^{s(1-\delta)}\Delta^{-\delta(1-s)}.
\end{equation*}

Putting this together returns the result.

\end{proof}

We can use this to pick out the optimal sequence of covering times.

\begin{cor}\label{Cor:N_n,i is best INT DIM}
    Let $\delta,s\in(0,1)$ and $\theta\in(0,1]$, let $n$ and $i$ satisfy the conditions of Lemma \ref{Lem:SingleLevelCost}, and let $N_{n,i}\leq N<N_{n,i+1}$. Then the infimum of $\sum_{k}|U_{k}|^{s}_{\dOm}$ over level-$N$ covers for $\theta$ is comparable to the infimum over level-$N_{n,i}$ covers for $\theta$, with constants depending on $\delta$. In particular, when bounding $\underline{\dim}_{\theta}\OmegaAD$ from above by single-level covers, it suffices to consider the times $N_{n,i}$.
\end{cor}
\begin{proof}
    By Lemma \ref{Lem:SingleLevelCost}, both infima are comparable to the same quantity, which depends on $n,i,\theta$ and $s$, but not on $N$ (with constants depending on just $\delta$). Since level-$N$ and level-$N_{n,i}$ covers for $\theta$ have the same scale, namely $\Delta_{n+2,i}$, any sequence of level-$N$ covers whose covering sums tend to $0$ may therefore be replaced by a sequence of level-$N_{n,i}$ covers for $\theta$ at the same scales whose covering sums also tend to $0$, for instance by taking the level-$N_{n,i}$ cover constructed in the proof of Lemma \ref{Lem:SingleLevelCost}, whose covering sum is at most a constant multiple of that of any level-$N$ cover.
 \end{proof}

\subsection{Upper bounds for \texorpdfstring{$\theta$}{theta}-intermediate dimension}\label{SubSec:UB}
 
\begin{thm}\label{Thm:lower Int Dim UB}
Let $\delta\in (0,1)$ and let $\alpha \in \Dioph(\nu) $ with $\nu>1$. Let $r_{\nu}$ be the unique positive root of
\begin{equation*}
    (1-\delta)r^2 +\delta r -\nu = 0.
\end{equation*}
Then, covering $\OmegaAD$ by a level-$N$ cover for $\theta$ with $N=N_{n,i}$ gives, for each $\theta\in (0,1]$, 
\begin{align*}
    \underline{\dim}_\theta\,\Omega^\delta_\alpha\;\le\;
\begin{cases}
\;\;\;\;\;\dfrac{\delta r_\nu}{\nu}, & 0<\theta\le\dfrac1{r_\nu},\\[3mm]
\dfrac{\delta\theta}{\delta\theta+1-\delta}, & \dfrac1{r_\nu}\le\theta\le 1,
\end{cases}
\end{align*}
and for each $\theta$ this estimate is attained by a level-$N$ cover for $\theta$ which minimises the covering sum, up to a constant depending only on $\delta$, among all level-$N$ covers for $\theta$ and all $N$.
\end{thm}

\begin{proof}    
    We begin by defining 
    \begin{equation*}
        \gamma_{n}:=\frac{\log q_{n+2}}{\log q_{n+1}}.
    \end{equation*}
    Then, since $\alpha\in\Dioph(\nu)$, there exists a subsequence along which $\gamma_{n}\rightarrow\nu$, and will therefore use that $q_{n+2}=q_{n+1}^{\gamma_{n}}$. In particular, since 
    $\nu>1$, we can pick a further subsequence such that $\gamma_{n}>1$ for all $n\in\mathbb{N}$. We pass to this subsequence and drop the subsequence notation. By Corollary \ref{Cor:N_n,i is best INT DIM}, it suffices to consider times $N=N_{n,i}$, and we will parametrise these times by writing
    \begin{equation*}
        N_{n,i}= \Theta\left( q_{n+1}^{r}\right), \;\; i=\Theta\left( q_{n+1}^{r-1}\right) 
    \end{equation*}
    which is legitimate since $N_{n,i}=\Theta\left( iq_{n+1}\right)$ and $i\in [1,a_{n+2}]$ with $a_{n+2}= \Theta\left( \frac{q_{n+2}}{q_{n+1}} \right)$. Here, $r=1$ corresponds to the two-gap time $N_{n}$ and $r=\gamma_{n}$ corresponds to the two-gap time $N_{n+1}$. For a fixed $r$, we choose $i=i(n)$ with $i(n)=\Theta(q_{n+1}^{r-1})$, and refer to the times $N_{n,i(n)}$, together with their covering scales $\Delta_{n+2,i(n)}$, as \emph{level $r$}.

    Let $r_{\gamma_{n}}$ denote the unique positive root of $(1-\delta)r^2+\delta r=\gamma_{n}$. Since the polynomial $(1-\delta)r^2+\delta r$ takes the value $1<\gamma_{n}$ at $r=1$  and a value exceeding $\gamma_{n}$ at $r=\gamma_{n}$, we have that $1<r_{\gamma_{n}}<\gamma_{n}$. In what follows, we will restrict to values of $r\in[1,r_{\gamma_{n}}]$, a range which contains the optimal level found below. So, uniformly in $r$ over this range, we have
    \begin{equation*}
        \frac{i}{a_{n+2}} = O\left( q_{n+1}^{r-\gamma_{n}}\right) \rightarrow0 \;\; \text{ for }r<\gamma_{n},
    \end{equation*}
    so we can further say that $i\leq \frac{a_{n+2}}{2}$ for all large $n$ in this subsequence, and we deduce that the conditions for Proposition \ref{Prop:BigAppendixProp} are met.

    We now consider the sizes of the complementary intervals in $\J_{N_{n,i}}$. By \eqref{Eqn:3.1}, we have that there are $q_{n+1}$ intervals of $q_{n+2,i}$-type, whose $\dOm$-length is 
    \begin{equation*}
        \Delta_{n+2,i}=\Theta\left( i^{1-1/\delta}q_{n+1}^{-1/\delta}\right)= \Theta \left(q_{n+1}^{-\frac{r(1-\delta)+\delta}{\delta}} \right),
    \end{equation*}
    and $q_{n}+iq_{n+1} = \Theta(q_{n+1}^r)$ intervals of $q_{n+1}$-type, whose $\dOm$-length is
    \begin{equation*}
        \Gamma_{n+1} = \Theta\left(q_{n+2}^{-1/\delta} \right)=\Theta\left(q_{n+1}^{-\frac{\gamma_{n}}{\delta}} \right).
    \end{equation*}
    Since $i\leq \frac{a_{n+2}}{2}$, we have that $\Gamma_{n+1}\leq \Delta_{n+2,i}$, so we fix our covering scale at $\Delta=\Delta_{n+2,i}:=\max\left\{|I|_{\dOm}:I \text{ is of } q_{n+2,i}\text{-type} \right\} $. The complementary intervals of $\J_{N_{n,i}}$ therefore form an admissible cover for $\theta$ when we have that $\Gamma_{n+1}
    \geq \Delta_{n+2,i}^{1/\theta}$, in other words, when 
    \begin{equation*}
        \frac{\gamma_{n}}{\delta}\leq \frac{r(1-\delta)+\delta}{\delta\theta},
    \end{equation*}
    or equivalently, when
    \begin{equation}\label{Eqn:ThetaBoundIntDimProof}
        \theta\le\frac{r(1-\delta)+\delta}{\gamma_{n}}.
    \end{equation}
    Suppose first that \eqref{Eqn:ThetaBoundIntDimProof} holds, and let each complementary interval of $\J_{N_{n,i}}$ be its own covering set $U_{j}$. Then
    \begin{align*}
        \sum_{j}|U_{j}|^s_{\dOm} = \Theta\left(q_{n+1}^{1-s\frac{r(1-\delta)+\delta}{\delta}} \right) + \Theta\left(q_{n+1}^{r-s\frac{\gamma_{n}}{\delta}} \right),
    \end{align*}
    which tends to $0$ exactly when both exponents are negative, in other words, when
    \begin{equation*}
        1-s\frac{r(1-\delta)+\delta}{\delta} <0 \;\; \text{ and } \;\; r-s\frac{\gamma_{n}}{\delta}<0,
    \end{equation*}
    which, upon solving for $s$, gives us that 
    \begin{equation}\label{Eqn:IntDimBoundsWhen3.2holds}
        \frac{\delta}{r(1-\delta)+\delta}<s \;\; \text{ and } \;\; \frac{\delta r}{\gamma_{n}}< s.
    \end{equation}
    
    Now we suppose that \eqref{Eqn:ThetaBoundIntDimProof} does not hold, and so $\Gamma_{n+1}<\Delta_{n+2,i}^{1/\theta}$. In this case, the $q_{n+1}$-type intervals are too small to be used as a covering set for our $\theta$-intermediate dimension estimates, so instead, we cover them with sets of size $\Delta_{n+2,i}^{1/\theta}$. The first term in the above covering sum is unchanged, but now, the second term becomes $\Theta\left(q_{n+1}^{r-s\frac{r(1-\delta)+\delta}{\delta\theta}} \right)$. In other words, we have
    \begin{equation*}
        \sum_{j}|U_{j}|^s_{\dOm} = \Theta\left(q_{n+1}^{1-s\frac{r(1-\delta)+\delta}{\delta}} \right) + \Theta\left(q_{n+1}^{r-s\frac{r(1-\delta)+\delta}{\delta\theta}} \right),
    \end{equation*}
    which tends to $0$ precisely when both exponents are negative, so when
    \begin{equation*}
        1-s\frac{r(1-\delta)+\delta}{\delta}<0 \;\; \text{ and } \;\; r-s\frac{r(1-\delta)+\delta}{\delta\theta}< 0
    \end{equation*}
    which gives us
    \begin{equation}\label{Eqn:IntDimBounds3.2NOTHold}
        \frac{\delta}{r(1-\delta)+\delta}<s \;\; \text{ and } \;\; \frac{\delta \theta r}{r(1-\delta)+\delta}< s.
    \end{equation}
    We note that the second condition in each of \eqref{Eqn:IntDimBoundsWhen3.2holds} and \eqref{Eqn:IntDimBounds3.2NOTHold} coincide exactly when we have equality in \eqref{Eqn:ThetaBoundIntDimProof}. This in turn means that, in both regimes, our covering sum will tend to $0$ along our chosen subsequence when we have $s$ greater than each of
    \begin{equation*}
        \frac{\delta}{r(1-\delta)+\delta}, \;\; \frac{\delta r}{\gamma_{n}}, \;\; \frac{\delta\theta r}{r(1-\delta)+\delta}.
    \end{equation*}
    In particular, when \eqref{Eqn:ThetaBoundIntDimProof} holds, we have that $\frac{\delta\theta r}{r(1-\delta)+\delta}\leq \frac{\delta r}{\gamma_{n}}$, and so $s>\frac{\delta r}{\gamma_{n}}$ implies that $s>\frac{\delta\theta r}{r(1-\delta)+\delta}$. Similarly, if \eqref{Eqn:ThetaBoundIntDimProof} does not hold, we have the opposite implication. In either case, we are only concerned with the interaction of one of these functions of $r$ and how it interacts with $\frac{\delta}{r(1-\delta)+\delta}$.
    We now want to determine the values of $r$ where we can minimise the maximum of the three functions above when estimating our lower $\theta$-intermediate dimension. We begin by noting that the first of the above quantities is strictly decreasing in $r$, while the latter two are strictly increasing in $r$ for all $r>0$. Since the maximum of the two increasing expressions is itself increasing in $r$, and the first expression is decreasing in $r$, the maximum of all three is minimised at the smaller of the two points where the first meets the other two. Equating the first two expressions gives us the equation
    \begin{equation*}
        \frac{\delta}{r(1-\delta)+\delta}=\frac{\delta r}{\gamma_{n}}
    \end{equation*}
     which can be rearranged to give
    \begin{equation*}
        (1-\delta)r^2+ \delta r=\gamma_{n}.
    \end{equation*}
    We recall that the unique positive root of this expression was denoted by $r_{\gamma_{n}}$. Equating the first and last expressions gives us
    \begin{equation*}
        \frac{\delta}{r(1-\delta)+\delta}=\frac{\delta\theta r}{r(1-\delta)+\delta} 
    \end{equation*}
    which rearranges to give
    \begin{equation*}
        r=\frac{1}{\theta}.
    \end{equation*}
    The optimal level is therefore found at $r=r^*$, where
    \begin{equation*}
        r^* :=\min\left\{r_{\gamma_{n}},\frac{1}{\theta} \right\}.
    \end{equation*}
  
    We note that $1< r_{\gamma_{n}}<\gamma_{n}$ since the polynomial $(1-\delta)r^2 +\delta r$ takes value $1$ at $r=1$, and a value exceeding $\gamma_{n}$ at $r=\gamma_{n}$. We also have that $1/\theta \geq 1$, so we must have that $1
    \leq r^* \leq \gamma_{n}$, meaning that $r^*$ is a legitimate candidate to be used in the dimension estimate. Hence, covering $\OmegaAD$ by complementary intervals of $\J_{N_{n,i}}$ with $r=r^{*}$, and letting $\gamma_{n}\rightarrow \nu$ as $n\rightarrow\infty$, we obtain
    \begin{equation*}
        \underline{\dim}_{\theta}\OmegaAD\leq \frac{\delta}{r^{*}(1-\delta)+\delta}.
    \end{equation*}
    If $\theta\leq \frac{1}{r_{\nu}}$, then $r^{*}=r_{\nu}$, so since $r_{\nu}(1-\delta)+\delta = \frac{\nu}{r_{\nu}}$, our $\theta$-intermediate dimension bound can be expressed as $\frac{\delta r_{\nu}}{\nu}$. If instead we have that $\theta\geq \frac{1}{r_{\nu}}$, then $r^*= \frac{1}{\theta}$, and the bound is 
    \begin{equation*}
        \frac{\delta}{\left(\frac{1}{\theta}\right)(1-\delta)+\delta} = \frac{\delta\theta}{\delta\theta+1-\delta}.
    \end{equation*}
Since $\frac{\delta\theta}{\delta\theta+1-\delta}$ is increasing in $\theta$ and equals $\frac{\delta r_{\nu}}{\nu}$ at $\theta=\frac{1}{r_{\nu}}$, in both cases this bound is the larger of $\frac{\delta r_{\nu}}{\nu}$ and $\frac{\delta\theta}{\delta\theta+1-\delta}$. 

It now remains to check that there are no level-$N$ covers for $\theta$ that do better than these estimates. By Corollary \ref{Cor:N_n,i is best INT DIM}, it suffices to consider level-$N_{n,i}$ covers for $\theta$, and by Lemma \ref{Lem:SingleLevelCost}, the infimal covering sum at such a level is of the form
\begin{equation*}
    \Theta\left(q_{n+1}\Delta^{s}_{n+2,i}+\min\left\{iq_{n+1}\beta^{s},\beta^{s(1-\delta)}\Delta_{n+2,i}^{-\delta(1-s)} \right\} \right),
\end{equation*}
where $\beta=\Gamma_{n+1}$ when \eqref{Eqn:ThetaBoundIntDimProof} holds and $\beta=\Delta_{n+2,i}^{1/\theta}$ otherwise.
The first two terms in this expression are exactly those of the covering sum computed above, and the argument shows that for every $r$, both terms tend to $0$ only if $s$ exceeds $\frac{\delta r_{\nu}}{\nu}$ and $\frac{\delta\theta}{\delta\theta+1-\delta}$. So it suffices to show that the same is true if the minimum in the statement of Lemma \ref{Lem:SingleLevelCost} takes the latter value. In particular, since each term in the covering sum expression can be given as a power of $q_{n+1}$, the infimal covering sum will tend to $0$ when both the first term tends to $0$, and at least one of the terms in the minimum tends to $0$. Each of these conditions is of the form $s>\cdot $, so the smallest $s$ for which the infimal covering sum at level $r$ tends to $0$ is 
\begin{multline}\label{Eqn:LevelThreshold}
    \max\Big\{ \text{threshold from the first term}, \\ \min\big\{ \text{threshold from the second term},\; \text{threshold from the third term} \big\} \Big\}.
\end{multline}
\pgfmathsetmacro{\dl}{0.5}
\pgfmathsetmacro{\gm}{3}
\pgfmathsetmacro{\rg}{2}
\pgfmathsetmacro{\crossval}{1/3}
 
\pgfplotsset{
  thresholdaxis/.style={
    width=0.62\textwidth,
    height=0.42\textwidth,
    clip=false,
    domain=1:\gm,
    samples=100,
    xmin=1, xmax=\gm,
    ymin=0, ymax=0.55,
    axis lines=left,
    xlabel={$r$},
    xtick={1,\rg,\gm},
    xticklabels={$1$,$r_{\gamma_{n}}$,$\gamma_{n}$},
    ytick={\crossval},
    yticklabels={$\frac{\delta r_{\gamma_{n}}}{\gamma_{n}}$},
  },
}
 
\begin{figure}[ht]
  \centering
  \begin{subfigure}[t]{\textwidth}
    \centering
    \begin{tikzpicture}
      \begin{axis}[thresholdaxis]
        \addplot[thick, solid] {\dl/(x*(1-\dl)+\dl)}
          node[pos=1, right, font=\small] {$\frac{\delta}{r(1-\delta)+\delta}$};
        \addplot[thick, dashed] {\dl*x/\gm}
          node[pos=1, right, font=\small] {$\frac{\delta r}{\gamma_{n}}$};
        \draw[dotted] (axis cs:\rg,0) -- (axis cs:\rg,\crossval);
        \draw[dotted] (axis cs:1,\crossval) -- (axis cs:\rg,\crossval);
        \fill (axis cs:\rg,\crossval) circle (1.5pt);
      \end{axis}
    \end{tikzpicture}
    \caption{The thresholds from the first and second terms.}
    \label{Fig:Threshold12}
  \end{subfigure}
  \par\bigskip
  \begin{subfigure}[t]{\textwidth}
    \centering
    \begin{tikzpicture}
      \begin{axis}[thresholdaxis]
        \addplot[thick, solid] {\dl/(x*(1-\dl)+\dl)}
          node[pos=1, right, font=\small] {$\frac{\delta}{r(1-\delta)+\delta}$};
        \addplot[thick, dash dot] {\dl*(x*(1-\dl)+\dl)/(\dl*(x*(1-\dl)+\dl)+(1-\dl)*\gm)}
          node[pos=1, right, font=\small] {$\frac{\delta(r(1-\delta)+\delta)}{\delta(r(1-\delta)+\delta)+(1-\delta)\gamma_{n}}$};
        \draw[dotted] (axis cs:\rg,0) -- (axis cs:\rg,\crossval);
        \draw[dotted] (axis cs:1,\crossval) -- (axis cs:\rg,\crossval);
        \fill (axis cs:\rg,\crossval) circle (1.5pt);
      \end{axis}
    \end{tikzpicture}
    \caption{The thresholds from the first and third terms, when \eqref{Eqn:ThetaBoundIntDimProof} holds.}
    \label{Fig:Threshold13}
  \end{subfigure}
  \caption{The thresholds as functions of $r$, for $\delta=\frac{1}{2}$ and $\gamma_{n}=3$, so that $r_{\gamma_{n}}=2$. In (B), $\theta\leq\frac{1}{3}$, so that \eqref{Eqn:ThetaBoundIntDimProof} holds for every $r\in[1,\gamma_{n}]$.}
  \label{Fig:Thresholds}
\end{figure}
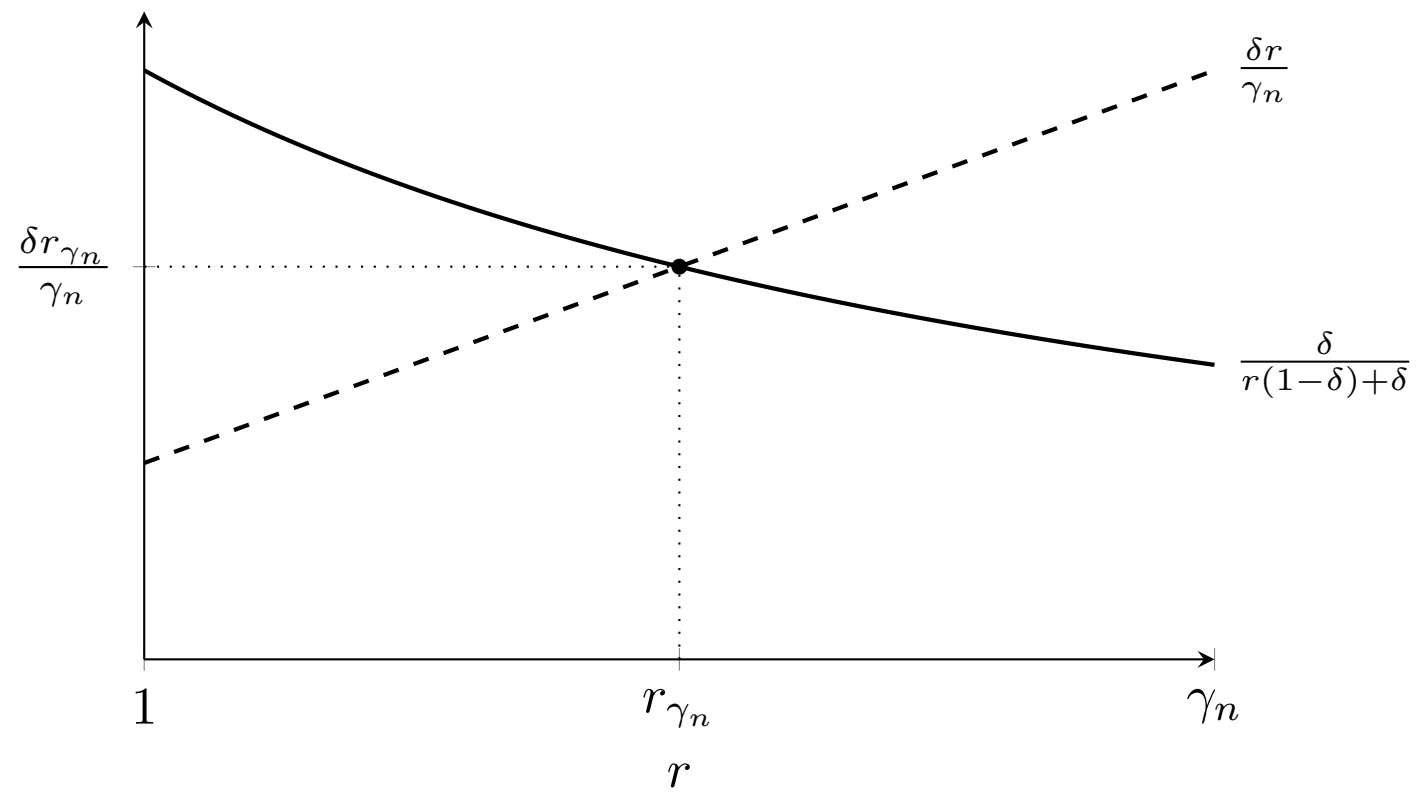
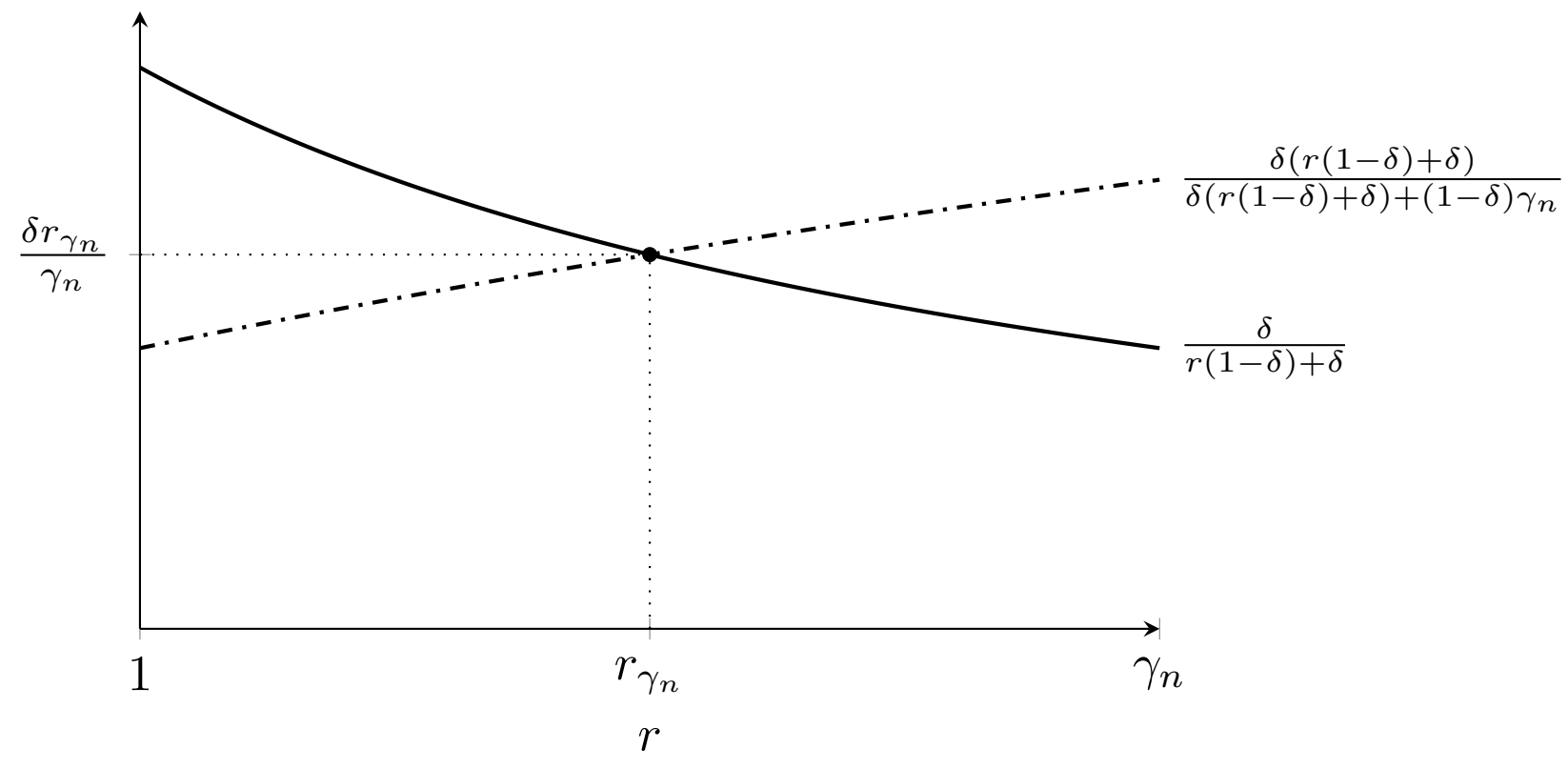
We have already computed the thresholds from the first two terms. The first term has threshold $\frac{\delta}{r(1-\delta)+\delta}$, which decreases in $r$, which, along with the threshold from the second term, gives the estimate stated in the theorem, namely the larger of $\frac{\delta r_{\gamma_{n}}}{\gamma_{n}}$ and $\frac{\delta\theta}{\delta\theta+1-\delta}$. We want to show that including the third term does not lower this smallest value.

We recall from earlier in the proof that 
\begin{equation*}
    \Gamma_{n+1}= \Theta\left(q_{n+1}^{-\frac{\gamma_{n}}{\delta}} \right) \;\; \text{ and } \;\; \Delta_{n+2,i}=\Theta\left( q_{n+1}^{-\frac{r(1-\delta)+\delta)}{\delta}}\right)
\end{equation*}
which means that when \eqref{Eqn:ThetaBoundIntDimProof} holds, we have that
\begin{equation*}
    \beta^{s(1-\delta)}\Delta_{n+2,i}^{-\delta(1-s)} = \Gamma_{n+1}^{s(1-\delta)}\Delta_{n+2,i}^{-\delta(1-s)} =\Theta\left( q_{n+1}^{(1-s)(r(1-\delta)+\delta) -s\frac{(1-\delta)\gamma_{n}}{\delta}}\right),
\end{equation*}
and when \eqref{Eqn:ThetaBoundIntDimProof} does not hold, we have that
\begin{equation*}
    \beta^{s(1-\delta)}\Delta_{n+2,i}^{-\delta(1-s)} = \Delta_{n+2,i}^{\frac{s(1-\delta)}{\theta}}\Delta_{n+2,i}^{-\delta(1-s)} = \Delta_{n+2,i}^{\frac{s(1-\delta)-\delta\theta(1-s)}{\theta}}.
\end{equation*}

We will consider each of these instances individually. If we are in the case where \eqref{Eqn:ThetaBoundIntDimProof} holds, we want to determine the value of $s$ for which
\begin{equation*}
    q_{n+1}^{(1-s)(r(1-\delta)+\delta) -s\frac{(1-\delta)\gamma_{n}}{\delta}} \rightarrow 0
\end{equation*}
as $n\rightarrow \infty$. This will happen when 
\begin{equation*}
    (1-s)(r(1-\delta)+\delta) -s\frac{(1-\delta)\gamma_{n}}{\delta} <0 \iff \frac{\delta(r(1-\delta)+\delta)}{\delta(r(1-\delta)+\delta)+(1-\delta)\gamma_{n}} <s.
\end{equation*}

We note that as a function of $r$, this is increasing. So, as before, we want to establish the smallest value of $r$ for which this meets the threshold from the first term; in other words, we want to find $r$ for which
\begin{equation*}
    \frac{\delta(r(1-\delta)+\delta)}{\delta(r(1-\delta)+\delta)+(1-\delta)\gamma_{n}}=\frac{\delta}{r(1-\delta)+\delta}
\end{equation*}
which holds exactly when
\begin{equation*}
    (r(1-\delta)+\delta)^{2}-\delta(r(1-\delta)+\delta)=(1-\delta)\gamma_{n}
\end{equation*}
which in turn is equivalent to
\begin{equation*}
    (r(1-\delta)+\delta)(r(1-\delta))=(1-\delta)\gamma_{n}.
\end{equation*}
Dividing through by $(1-\delta)$ gives us
\begin{equation*}
    r(r(1-\delta)+\delta)=\gamma_{n}
\end{equation*}
which is exactly the polynomial for which $r_{\gamma_{n}}$ is the unique positive solution. So, the increasing thresholds from both the second and third terms meet the decreasing threshold from the first at $r=r_{\gamma_{n}}$ at height $\frac{\delta r_{\gamma_{n}}}{\gamma_{n}}$. We note that when $r>r_{\gamma_{n}}$, the thresholds from the second and third terms, which are increasing, force $s>\frac{\delta r_{\gamma_{n}}}{\gamma_{n}}$, and similarly, for $r<r_{\gamma_{n}}$ we have $s>\frac{\delta}{r(1-\delta)+\delta}>\frac{\delta r_{\gamma_{n}}}{\gamma_{n}}$. In other words, sharing covering sets may be cheaper to the covering sum than using individual covers at levels $r>r_{\gamma_{n}}$, but it does not change the level $r_{\gamma_{n}}$ at which the smallest possible dimension estimate is extracted. Finally, we note that when \eqref{Eqn:ThetaBoundIntDimProof} holds, we have $r(1-\delta)+\delta\geq \theta\gamma_{n}$, so since the threshold for the third term is increasing in $r(1-\delta)+\delta$, we have that when \eqref{Eqn:ThetaBoundIntDimProof} holds,
\begin{equation*}
    \frac{\delta(r(1-\delta)+\delta)}{\delta(r(1-\delta)+\delta)+(1-\delta)\gamma_{n}} \geq \frac{\delta\theta\gamma_{n}}{\delta\theta\gamma_{n}+(1-\delta)\gamma_{n}} = \frac{\delta\theta}{\delta\theta+1-\delta}.
\end{equation*}
In either case, we have shown that in this regime, the dimension estimate extracted from \eqref{Eqn:LevelThreshold} is at least the larger of $\frac{\delta r_{\gamma_{n}}}{\gamma_{n}}$ and $\frac{\delta\theta}{\delta\theta+1-\delta}$ at every level.

Suppose \eqref{Eqn:ThetaBoundIntDimProof} does not hold. Then $\beta=\Delta^{1/\theta}_{n+2,i}$, and since $\Delta_{n+2,i}\rightarrow 0$ as $n\rightarrow\infty$, the third term in the covering sum will tend to zero when the exponent is positive, in other words, when
\begin{equation*}
    \frac{s(1-\delta)-\delta\theta(1-s)}{\theta}>0 \iff s>\frac{\delta\theta}{\delta\theta+1-\delta}.
\end{equation*}
This threshold is independent of $r$, and agrees with the second branch estimate. It remains to account for the first term. Since \eqref{Eqn:ThetaBoundIntDimProof} does not hold, we have that $\theta\gamma_{n}>r(1-\delta)+\delta$. In particular, this restricts the level at which we can view the system, constraining $r$. The first threshold then requires
\begin{equation*}
    s>\frac{\delta}{r(1-\delta)+\delta}>\frac{\delta}{\theta\gamma_{n}}.
\end{equation*}
If $\theta\geq\frac{1}{r_{\gamma_{n}}}$, then the dimension estimate from the first two terms is at best $\frac{\delta\theta}{\delta\theta+1-\delta}$, and so introducing the third term does not improve this. If instead we have $\theta< \frac{1}{r_{\gamma_{n}}}$, then $\theta\gamma_{n}<\frac{\gamma_{n}}{r_{\gamma_{n}}}=r_{\gamma_{n}}(1-\delta)+\delta$, so the first term requires $s>\frac{\delta r_{\gamma_{n}}}{\gamma_{n}}$. So the third term permits every $s$ above the previous estimate. So in this regime too, \eqref{Eqn:LevelThreshold} is at least $\max\left\{\frac{\delta r_{\gamma_{n}}}{\gamma_{n}} ,\frac{\delta\theta}{\delta\theta+1-\delta}\right\}$ at every level.

Hence, at each level, the infimal covering sum tends to $0$ only if $s>\max\left\{\frac{\delta r_{\gamma_{n}}}{\gamma_{n}} ,\frac{\delta\theta}{\delta\theta+1-\delta}\right\}$. Since each exponent of $q_{n+1}$ computed above is an affine function of $s$ which vanishes at the corresponding threshold and has slope at most $-1$, this conclusion is uniform in $r$. I.e. for $s$ below this maximum, the infimal covering sum at every level is at least 
\begin{equation*}
    cq_{n+1}^{\max\left\{\frac{\delta r_{\gamma_{n}}}{\gamma_{n}} ,\frac{\delta\theta}{\delta\theta+1-\delta}\right\}-s}
\end{equation*}
where $c>0$ does not depend on $r$ or $n$. The argument only depends on $n$ through $\gamma_{n}$, and by Lemma \ref{lem:DiophClass} we have $\gamma_{n}\leq \nu+\eps$ for all sufficiently large $n$, while $\frac{\delta r_{\gamma}}{\gamma}=\frac{\delta}{r_{\gamma}(1-\delta)+\delta}$ is continuous and decreasing in $\gamma$. Therefore, for every $s<\max\left\{\frac{\delta r_{\gamma_{n}}}{\gamma_{n}} ,\frac{\delta\theta}{\delta\theta+1-\delta}\right\}$, the infimal covering sum over level-$N$ covers for $\theta$ tends to infinity as the covering scales tend to $0$, and no level-$N$ cover for $\theta$ improves on the estimate above.

\end{proof}

\begin{cor}\label{Cor:dimH UB}
    Let $\delta\in(0,1)$ and $\alpha\in\Dioph(\nu)$ with $\nu>1$. Then $\dim_{\mathrm{H}}\OmegaAD\leq \frac{\delta r_{\nu}}{\nu}$. Moreover, $1<r_{\nu}<\nu$, so this value lies between $\frac{\delta}{\nu}$ and $\delta$.
\end{cor}
\begin{proof}
    Since $\dim_{\mathrm{H}}\OmegaAD\leq \underline{\dim}_{\theta}\OmegaAD$ for every $\theta\in(0,1]$, the bound follows from Theorem \ref{Thm:lower Int Dim UB} with any $\theta\leq \frac{1}{r_{\nu}}$. The fact that $1<r_{\nu}<\nu$ follows as in the proof of Theorem \ref{Thm:lower Int Dim UB}, with $\nu$ in place of $\gamma_{n}$.
\end{proof}

\subsection{Lower bounds for lower \texorpdfstring{$\theta$}{theta}-intermediate dimension}\label{SubSec:LB}
We now turn our attention to finding a lower bound for $\underline{\dim}_{\theta}\OmegaAD$. We make use of a result found in \cite[Corollary 3.14]{Banaji2023} relating the $\theta$-intermediate dimension to the \emph{Assouad dimension}. Similar results can be found in \cite{Falconer2021}. First, we must define the Assouad dimension. For a thorough study of the Assouad dimension, we refer the reader to \cite{fraser2020assouaddimensionfractalgeometry}.

\begin{definition}
    Suppose a subset $F$ of a metric space has more than one point. Then the \emph{Assouad dimension} of $F$ is defined by

\begin{align*}
        \dim_{\mathrm{A}}F:=\inf\bigg\{a:  &\text{  there exists }C>0 \text{ such that } N_{r}(B_{R}(x)\cap F)\leq C\left(\frac{R}{r}\right)^{a}  \\   &\;\;\;\;\;\;\text{ for all } x\in F \text{ and } 0<r<R \bigg\}.
    \end{align*}
\end{definition}

The Assouad dimension gives information about the `thickest' part of a set. We now state the following proposition.

\begin{prop}\cite[Corollary 3.14]{Banaji2023}
    For a bounded set $F\subset \mathbb{R}^n$ and $\theta\in(0,1]$, 
    \begin{equation*}
        \underline{\dim}_{\theta}F\geq \frac{\theta\dim_{\mathrm{A}}F\,\underline{\dim}_{\mathrm{B}}F}{\dim_{\mathrm{A}}F-(1-\theta)\underline{\dim}_{\mathrm{B}}F}
    \end{equation*}
    with a similar inequality for upper dimensions. 
\end{prop}

The right-hand side of the bound above is decreasing in $\dim_{\mathrm{A}}F$, so it remains a lower bound when $\dim_{\mathrm{A}}F$ is replaced by any upper bound for it. This gives the following.

\begin{cor}
    For a bounded set $F\subset \mathbb{R}^n$ with $\dim_{\mathrm{A}}F\leq 1$ and $\theta\in(0,1]$, 
    \begin{equation}\label{Eqn:IntDimLB Banaji}
        \underline{\dim}_{\theta}F\geq \frac{\theta\,\underline{\dim}_{\mathrm{B}}F}{1-(1-\theta)\underline{\dim}_{\mathrm{B}}F}
    \end{equation}
    with a similar inequality for upper dimensions. 
\end{cor}
The arc-length metric on $S^1$ is bi-Lipschitz equivalent to the Euclidean metric on a circle in $\mathbb{R}^2$, which has Assouad dimension at most $1$.
With this in hand, we immediately get a lower bound on the lower $\theta$-intermediate dimension of $\OmegaAD$.

\begin{thm}\label{Thm:IntDim LB}
    Let $\delta\in(0,1)$ and $\alpha\in \Dioph(\nu)$ for $\nu>1$. Then for $\theta\in (0,1]$, 
    \begin{equation*}
        \underline{\dim}_{\theta}\OmegaAD\geq \max\left\{\dim_{\mathrm{H}}\OmegaAD,\frac{\delta\theta}{\delta\theta+1-\delta}\right\}, 
    \end{equation*}
    and the second term is larger whenever $\theta\geq \frac{1}{r_{\nu}}$.
\end{thm}
\begin{proof}
    From \cite{KraSchmeling}, we have that $\dim_{\mathrm{B}}\OmegaAD=\underline{\dim}_{\mathrm{B}}\OmegaAD=\delta$, then, substituting this into \eqref{Eqn:IntDimLB Banaji}, gives us that $\underline{\dim}_{\theta}\OmegaAD\geq\frac{\delta\theta}{\delta\theta+1-\delta}$. We note that $\underline{\dim}_{\theta}\OmegaAD\geq \dim_{\mathrm{H}}\OmegaAD$, and for $\theta\geq\frac{1}{r_{\nu}}$, we have that $\frac{\delta\theta}{\delta\theta+1-\delta}\geq \frac{\delta r_{\nu}}{\nu}\geq \dim_{\mathrm{H}}\OmegaAD$, by Corollary \ref{Cor:dimH UB}. 
\end{proof}

Theorems \ref{Thm:lower Int Dim UB} and \ref{Thm:IntDim LB} together give Theorem \ref{Thm:main}. We do not know whether either bound in Theorem \ref{Thm:main} for $\theta<1/r_{\nu}$ is sharp. We discuss what we do know below.

\begin{rmk}
    For $\theta\in(0,1/r_{\nu})$, Theorems \ref{thm:KS}, \ref{Thm:lower Int Dim UB} and \ref{Thm:IntDim LB} tell us that 
    \begin{equation*}
        \max \left\{ \frac{\delta}{\nu}, \frac{\delta\theta}{\delta\theta+1-\delta}\right\} \leq \underline{\dim}_{\theta}\OmegaAD\leq\frac{\delta r_{\nu}}{\nu}.
    \end{equation*}
    Since $\frac{\delta\theta}{\delta\theta+1-\delta}=\frac{\delta r_{\nu}}{\nu}$ when $\theta=1/r_{\nu}$, the two bounds agree at this point, but differ for all $\theta<1/r_{\nu}$ since $\frac{\delta\theta}{\delta\theta+1-\delta}$ decreases with $\theta$, and we do not know whether either bound is sharp.
    
    We can make some observations. Since $\underline{\dim}_{\theta}\OmegaAD$ is non-decreasing in $\theta$, if the upper bound is attained at some $\theta_{0}\in (0,1/r_{\nu})$, then it is attained for all $\theta\in[\theta_{0},1/r_{\nu})$. In particular, the upper bound is sharp on the whole range of $\theta$ if and only if $\lim_{\theta\rightarrow 0} \underline{\dim}_{\theta}\OmegaAD=\frac{\delta r_{\nu}}{\nu}$.

    Kra and Schmeling state in \cite[Corollary 3.2]{KraSchmeling} that $\dim_{\mathrm{H}}\OmegaAD=\frac{\delta}{\nu}$. If this holds, then either $\underline{\dim}_{\theta}\OmegaAD<\frac{\delta r_{\nu}}{\nu}$ for all sufficiently small $\theta$, or $\underline{\dim}_{\theta}\OmegaAD$ is discontinuous at $0$.

    We showed in Theorem \ref{Thm:lower Int Dim UB} that no cover of complementary intervals of $\J_{N}$ at a single level $N$ (as defined in Definition \ref{Def:LevelNCover}) can improve on $\frac{\delta r_{\nu}}{\nu}$. To improve on this, we therefore would have to consider covers across multiple levels. An example of how this might look would be to cover the long intervals of $\J_{N_{n,i}}$ at some level $N_{n,i}$, then cover the clustered smaller intervals using complementary intervals of $\J_{N'}$ for some $N'>N_{n,i}$. This is similar to Kra and Schmeling's approach for their Hausdorff dimension upper bound in \cite{KraSchmeling}, where they separate regions which experience different distortion under the semiconjugacy $h$. For the intermediate dimension, the restriction $\Delta^{1/\theta}\leq|U_{k}|_{\dOm}\leq \Delta$ limits how far these refinements can go, so quantifying this limitation should equate to showing whether this upper bound is sharp. To refine the lower bound, a natural approach would be to use the mass distribution principle for intermediate dimensions \cite{Falconer2018}.

\end{rmk}

We end with a remark on the upper intermediate dimensions.

\begin{rmk}\label{Rmk:UpperLower}
Since $\overline{\dim}_{\theta}\geq\underline{\dim}_{\theta}$, Theorem \ref{Thm:main} gives
\begin{equation*}
    \frac{\delta\theta}{\delta\theta+1-\delta}\leq\overline{\dim}_{\theta}\OmegaAD\leq\delta, \qquad \theta\in(0,1],
\end{equation*}
where the upper bound holds as $\overline{\dim}_{\theta}\OmegaAD\leq\overline{\dim}_{\mathrm{B}}\OmegaAD=\delta$.  Angelini and Molter \cite{AngeliniMolter} show that the lower bound depends only on the lengths of the gaps: it is the least possible upper $\theta$-intermediate dimension of a closed set whose complementary intervals have lengths $(l_{n})_{n\in\mathbb{Z}}$, attained by the countable set of tail sums of these lengths. They also show that, for a broad class of gap sequences, every value between this minimum and the box dimension is attained by some arrangement of the gaps. Further, in \cite{Angelini2026UpperIntDimComp}, Angelini characterises the attainable forms of the upper intermediate dimension of complementary sets, linking the result to the work of Banaji and Rutar in \cite{BanajiRutar}. Determining where $\overline{\dim}_{\theta}\OmegaAD$ lies in this range requires efficient covers at every small scale, rather than along a subsequence as in the proof of Theorem \ref{Thm:main}, and we do not pursue this here.
\end{rmk}

\bibliography{biblio}
\bibliographystyle{alpha}

\end{document}